\documentclass[12pt,reqno]{amsart}
\usepackage[T1]{fontenc}
\usepackage{lmodern}
\usepackage{amsmath,amssymb,mathtools,amsthm, amscd, amsfonts, graphicx, color,mathrsfs,enumerate}
\usepackage[margin=1in]{geometry}
\usepackage[bookmarksnumbered, colorlinks, plainpages,linkcolor=black,anchorcolor=black,citecolor=black,urlcolor=black]{hyperref}
\usepackage{exscale}
\usepackage[all]{xy}
\usepackage{relsize}
\usepackage{pdfsync}

\makeatletter
\@namedef{subjclassname@2020}{%
	\textup{2020} Mathematics Subject Classification}
\makeatother

\newtheorem{theorem}{Theorem}[section]
\newtheorem{lemma}[theorem]{Lemma}

\theoremstyle{definition}

\newtheorem{example}[theorem]{Example}

\theoremstyle{remark}

\numberwithin{equation}{section}

\newcommand{\B}{\mathbb B_n}
\newcommand{\bloch}{\mathcal B}
\newcommand{\Sph}{\mathbb S_n}
\newcommand{\D}{\mathbb D}
\newcommand{\T}{\mathbb T}
\newcommand{\BT}{\mathcal{BT}}
\newcommand{\HT}{\mathcal{HT}}
\newcommand{\hol}{\mathcal H}
\newcommand{\M}{\mathcal{M}}
\newcommand{\N}{\mathcal{N}}
\newcommand{\Q}{\mathcal{Q}}
\newcommand{\dd}{\,d}
\newcommand{\ip}[2]{\langle #1,#2\rangle}
\newcommand{\norm}[1]{\left\lVert #1\right\rVert}

\newcommand{\less}{\lesssim}

\DeclareMathOperator*{\esssup}{ess\,sup}

\begin{document}
\title[Multipliers and integration operators]{Multipliers and integration operators on Hardy--Bloch type tent spaces}
\date{November 21, 2026}
\author[J. Chen]{Jiale Chen}
\address{Jiale Chen, School of Mathematics and Statistics, Shaanxi Normal University, Xi'an 710119, China.}
\email{jialechen@snnu.edu.cn}

\thanks{This work was supported by National Natural Science Foundation (No. 12501170), Young Talent Fund of Association for Science and Technology in Shaanxi (No. 20260506) and the Fundamental Research Funds for the Central Universities (No. GK202601003) of China.}

\begin{abstract}
We investigate the boundedness of integration operators $J_g$ and $I_g$ and pointwise multipliers $M_g$ on Hardy--Bloch type tent spaces $\mathcal{BT}^p$. The main difficulty is that these spaces are defined through a non-tangential maximal function of the radial derivative, which places the problem at the interface of holomorphic function theory and harmonic analysis. We completely characterize the boundedness of $I_g:\mathcal{BT}^p\to\mathcal{BT}^q$ for all $0<p,q<\infty$, and obtain a complete characterization for $J_g$ and $M_g$ when $p<q$. For $p\geq q$, we establish separate necessary and sufficient conditions. The proofs in the case $p\leq q$ are mainly based on some fractional maximal function estimates, whereas those in the case $p>q$ are mainly based on a synthesis lemma. We also give examples showing that several converse implications in the range $p\geq q$ fail.
\end{abstract}
\maketitle


\section{Introduction}
\allowdisplaybreaks[4]

Let $\B$ be the unit ball of $\mathbb C^n$, and let $\Sph:=\partial\B$ be the unit sphere. For $\gamma>1$ and $\xi\in\Sph$, let $\Gamma_{\gamma}(\xi)$ be the admissible Kor\'anyi approach region given by
$$\Gamma_{\gamma}(\xi):=
\left\{z\in\B:|1-\ip{z}{\xi}|<\frac{\gamma}{2}(1-|z|^2)\right\}.$$
Throughout the paper, we write $\Gamma(\xi):=\Gamma_{2}(\xi)$. Let $0<p,q<\infty$ and $\alpha\in\mathbb{R}$. The tent space $T_{q,\alpha}^p$ consists of measurable functions $f$ on $\B$ such that
$$\|f\|^p_{T_{q,\alpha}^p}:=
\int_{\Sph}\left(\int_{\Gamma(\xi)}|f(z)|^q(1-|z|^2)^{\alpha}\dd V(z)\right)^{p/q}\dd\sigma(\xi)<\infty,$$
where $dV$ and $d\sigma$ are the normalized Lebesgue measure on $\B$ and the normalized surface measure on $\Sph$, respectively. In addition, the tent space $T^p_{\infty}$ consists of measurable functions $f$ on $\B$ with
$$\|f\|_{T^p_{\infty}}^p:=\int_{\Sph}\left(\esssup_{z\in\Gamma(\xi)}
|f(z)|\right)^p\dd\sigma(\xi)<\infty,$$
and the tent space $T^{\infty}_{q,\alpha}$ consists of measurable functions $f$ on $\B$ such that
$$\|f\|_{T^{\infty}_{q,\alpha}}:=\esssup_{\xi\in\Sph}
\left(\sup_{u\in\Gamma(\xi)}\frac{1}{(1-|u|^2)^n}
\int_{Q(u)}|f(z)|^q(1-|z|^2)^{n+\alpha}\dd V(z)\right)^{1/q}<\infty,$$
where $Q(0)=\B$ and for $u\in\B\setminus\{0\}$,
$$Q(u):=\{z\in\B:|1-\ip{z}{u/|u|}|<1-|u|^2\}.$$
In the above definitions, the Kor\'anyi region $\Gamma(\xi)$ can be replaced by $\Gamma_{\gamma}(\xi)$ for any aperture $\gamma>2$, and it is well-known that all fixed apertures $\gamma\geq2$ yield the same function space with equivalent (quasi-)norms. These tent spaces were first introduced by Coifman, Meyer and Stein \cite{CMS} in order to study several problems in harmonic analysis. Since then, these spaces have become powerful tools in harmonic analysis and operator theory; see for instance \cite{CN,Lu,LP-2,LPW,MPPW,PP,PR}.

In 2018, Per\"al\"a \cite{Pe} initiated the investigation of holomorphic tent spaces on the unit ball $\B$. Let $\hol(\B)$ be the space of holomorphic functions on $\B$. For $0<p,q<\infty$ and $\alpha>-n-1$, the holomorphic Hardy type tent space $\HT_{q,\alpha}^p$ and the Hardy--Carleson type tent space $\mathcal{CT}_{q,\alpha}$ are defined respectively by
$$\HT_{q,\alpha}^p:=T_{q,\alpha}^p\cap\hol(\B)$$
and
$$\mathcal{CT}_{q,\alpha}:=T_{q,\alpha}^{\infty}\cap\hol(\B)$$
with the corresponding (quasi-)norms inherited. The Hardy--Bloch type tent space $\BT^p$ consists of holomorphic functions $f$ on $\B$ such that $Rf(\cdot)(1-|\cdot|^2)\in T_{\infty}^p$, with the (quasi-)norm
$$\|f\|_{\BT^p}:=|f(0)|+\big\|Rf(\cdot)(1-|\cdot|^2)\big\|_{T_{\infty}^p}.$$
Here and in the sequel, for $f\in\hol(\B)$, $Rf$ denotes its radial derivative, given by
$$Rf(z):=\sum_{k=1}^nz_k\frac{\partial f}{\partial z_k}(z),\quad
z=(z_1,\dots,z_n)\in\B.$$
It is quite interesting that these holomorphic tent spaces include several classical spaces of holomorphic functions on the unit ball $\B$. For example, consider the Hardy space $H^p$ consisting of $f\in\hol(\B)$ such that
$$\|f\|_{H^p}^p:=\sup_{0<r<1}\int_{\Sph}|f(r\xi)|^p\dd\sigma(\xi)<\infty.$$
Then the Calder\'on area theorem (see \cite{Ca,FS} or \cite[Theorem B]{Pa}) indicates that $f\in H^p$ if and only if $Rf\in\HT^p_{2,1-n}$. Moreover, it was pointed out in \cite{Pe} (see also \cite[Proposition 2.4]{CN}) that the continuous embedding $H^p\subset\BT^p$ holds. Per\"al\"a \cite{Pe} investigated the dual of Hardy type tent spaces $\HT_{q,\alpha}^p$ for all $0<p,q<\infty$ and $\alpha>-n-1$, where Hardy--Carleson type tent spaces $\mathcal{CT}_{q,\alpha}$ and Hardy--Bloch type tent spaces $\BT^p$ emerged naturally. Subsequently, operator theory on the spaces $\HT_{q,\alpha}^p$ and $\mathcal{CT}_{q,\alpha}$ developed quickly; see \cite{AG,AMPR,Ch24,CW22,WZ21,WZ22}. However, to the best of our knowledge, there appear to be no operator-theoretic results on the spaces $\BT^p$ so far. In this paper, we are going to study multipliers and integration operators on the spaces $\BT^p$.

Given $g\in\hol(\B)$, the Volterra type integration operators $J_g$ and $I_g$ are respectively defined by
$$J_gf(z):=\int_0^1 f(tz)Rg(tz)\,\frac{dt}{t},\quad z\in\B,\quad f\in\hol(\B),$$
and
$$I_gf(z):=\int_0^1 Rf(tz)g(tz)\,\frac{dt}{t},\quad z\in\B,\quad f\in\hol(\B).$$
We also use $M_g$ to denote the pointwise multiplier, defined by
$$M_gf:=gf,\quad f\in\hol(\B).$$
Then $R(J_gf)=fRg$, $R(I_gf)=gRf$, and integration by parts gives
\begin{equation}\label{eq:identities}
M_gf=f(0)g(0)+J_gf+I_gf.
\end{equation}
These operators, called analytic paraproducts in \cite{ACFPP} due to the above relation, have been studied extensively on various spaces of holomorphic functions; see  \cite{CPW,CW22,GGP,GP,MPPW,Pa,WZ22} and the references therein. In particular, Wang and Zhou \cite{WZ22} established some embedding theorems for Hardy type tent spaces $\HT^p_{q,\alpha}$ and characterized the boundedness of the integration operators $J_g$ on these spaces. Chen and Wang \cite{CW22} studied the boundedness, compactness and strict singularity of $J_g$ acting on Hardy--Carleson type tent spaces $\mathcal{CT}_{q,\alpha}$.

The purpose of this paper is to investigate the boundedness of $J_g$, $I_g$ and $M_g$ acting on Hardy--Bloch type tent spaces $\BT^p$. The main difficulty is that the spaces $\BT^p$ are defined through a non-tangential maximal function of the radial derivative, rather than through the function itself or an area-integral functional. Thus their geometry is closely connected with harmonic analysis, including maximal operators, tent regions, and Carleson-type estimates. At the same time, the space $\BT^p$ combines features of Hardy and Bergman spaces by coupling boundary $L^p$ control with a Bloch-type derivative condition in the interior. Consequently, arguments from either Hardy spaces or Bergman spaces alone do not directly apply, and the proofs require a combination of their ideas together with estimates adapted to the derivative maximal function.

To state our main results, we introduce some function spaces. For $\beta\in\mathbb{R}$, the Bloch type space $\bloch_{\beta}$ consists of $f\in\hol(\B)$ such that
$$\|f\|_{\bloch_{\beta}}:=|f(0)|+\sup_{z\in\B}|Rf(z)|(1-|z|^2)^{\beta}<\infty.$$
When $\beta=1$, we obtain the classical Bloch space $\bloch:=\bloch_1$. It is clear that $\bloch_{\beta}=\mathbb{C}$ if $\beta<0$. The logarithmic Bloch space $\bloch_{\log}$ consists of $f\in\hol(\B)$ with
$$\|f\|_{\bloch_{\log}}:=|f(0)|+\sup_{z\in\B}|Rf(z)|(1-|z|^2)\log\frac{2}{1-|z|^2}<\infty.$$
For $0<p<\infty$, let $\BT^p_{\log}$ be the logarithmically weighted Hardy--Bloch type tent space, consisting of $f\in\hol(\B)$ such that
$$\|f\|_{\BT^p_{\log}}:=|f(0)|+\left(\int_{\Sph}\left(
\sup_{z\in\Gamma(\xi)}|Rf(z)|(1-|z|^2)
\log\frac{2}{1-|z|^2}\right)^p\dd\sigma(\xi)\right)^{1/p}<\infty.$$
We also use $H^{\infty}$ to denote the space of bounded holomorphic functions on $\B$ equipped with the supremum norm.

We are now ready to state our main results. In the case $p<q$, we completely characterize the boundedness of $J_g$, $I_g$ and $M_g$ acting from $\BT^p$ into $\BT^q$.

\begin{theorem}\label{main1}
Let $0<p<q<\infty$ and $g\in\hol(\B)$. Then:
\begin{enumerate}
	\item [(1)] $I_g:\BT^p\to\BT^q$ is bounded if and only if $g\equiv0$.
	\item [(2)] $M_g:\BT^p\to\BT^q$ is bounded if and only if $g\equiv0$.
	\item [(3)] $J_g:\BT^p\to\BT^q$ is bounded if and only if $g\in\bloch_{1+n/q-n/p}$.
\end{enumerate}
\end{theorem}

In the case $p\geq q$, the boundedness of $I_g$ acting from $\BT^p$ into $\BT^q$ is completely characterized, and separate necessary and sufficient conditions are established for the boundedness of $J_g$ and $M_g$.

\begin{theorem}\label{main2}
Let $0<p<\infty$ and $g\in\hol(\B)$. Then:
\begin{enumerate}
	\item [(1)] $I_g:\BT^p\to\BT^p$ is bounded if and only if $g\in H^{\infty}$.
	\item [(2)] If $J_g:\BT^p\to\BT^p$ is bounded, then $g\in\bloch$.
	\item [(3)] If $g\in\bloch_{\log}$, then $J_g:\BT^p\to\BT^p$ is bounded.
	\item [(4)] If $M_g:\BT^p\to\BT^p$ is bounded, then $g\in H^{\infty}$.
	\item [(5)] If $g\in H^{\infty}\cap\bloch_{\log}$, then $M_g:\BT^p\to\BT^p$ is bounded.
\end{enumerate}
\end{theorem}

\begin{theorem}\label{main3}
Let $0<q<p<\infty$ and $g\in\hol(\B)$. Then:
\begin{enumerate}
	\item [(1)] $I_g:\BT^p\to\BT^q$ is bounded if and only if $g\in H^{\frac{pq}{p-q}}$.
	\item [(2)] If $J_g:\BT^p\to\BT^q$ is bounded, then $g\in\BT^{\frac{pq}{p-q}}$.
	\item [(3)] If $g\in\BT^{\frac{pq}{p-q}}_{\log}$, then $J_g:\BT^p\to\BT^q$ is bounded.
	\item [(4)] If $M_g:\BT^p\to\BT^q$ is bounded, then
	$g\in H^{\frac{pq}{p-q}}$.
	\item [(5)] If $g\in H^{\frac{pq}{p-q}}\cap\BT_{\log}^{\frac{pq}{p-q}}$, then
	$M_g:\BT^p\to\BT^q$ is bounded.
	\item [(6)] The following two conditions are equivalent:
	\begin{enumerate}
		\item [(a)] $M_g:\BT^p\to\BT^q$ is bounded;
		\item [(b)] $g\in H^{\frac{pq}{p-q}}$ and $J_g:\BT^p\to\BT^q$ is bounded.
	\end{enumerate}
\end{enumerate}
\end{theorem}

The proofs use several complementary devices to overcome the preceding difficulty. In the case $p\leq q$, we establish two embedding theorems for $\BT^p$ (see Theorems \ref{embedding} and \ref{logembedding}) via some fractional maximal function estimates, which lead to the sufficiency parts in Theorems \ref{main1} and \ref{main2}, whereas the necessity parts follow from some pointwise estimates for functions in $\BT^p$ and a careful analysis of the geometry of Kor\'anyi regions $\Gamma(\xi)$. In the case $p>q$, the proof is mainly based on a synthesis lemma (Lemma \ref{synthesis}) and an approximation argument. The behavior of the operators $J_g$ and $M_g$ in the case $p\geq q$ is substantially more delicate than in the case $p<q$. In particular, we will construct some examples via lacunary series and Blaschke products on the unit disk $\D=\mathbb{B}_1$ in Section \ref{example} to show that, for the nonequivalent assertions in Theorems \ref{main2} and \ref{main3}, the converse implications are in general false.

This paper is organized as follows. In Section \ref{elementary}, we give some pointwise estimates for functions in $\BT^p$, and then establish some necessary conditions for the boundedness of multipliers and integration operators. Section \ref{proofs of 12} is devoted to the proofs of Theorems \ref{main1} and \ref{main2}, and Theorem \ref{main3} is proved in Section \ref{proof of 3}. Finally, we present some examples in Section \ref{example}.

Throughout the paper, we use $A\less B$ (or $B\gtrsim A$) to denote that there exists a positive constant $C$ independent of the functions and points under consideration such that $A\leq CB$. If $A\less B\less A$, then we write $A\asymp B$.

\section{Pointwise estimates and elementary necessary conditions}\label{elementary}

In this section, we first give some growth estimates for functions in Hardy--Bloch type tent spaces, and then establish some necessary conditions for the boundedness of multipliers and integration operators. We begin with the following notation. For $f\in\hol(\B)$ and $\gamma>1$, define
$$\Q_{\gamma}f(\xi):=
\sup_{z\in\Gamma_{\gamma}(\xi)}|Rf(z)|(1-|z|^2),\quad \xi\in\Sph.$$
When $\gamma=2$, we simply write $\Q f:=\Q_2f$. By the definition of the spaces $\BT^p$, a function $f\in\hol(\B)$ belongs to $\BT^p$ if and only if $\Q f\in L^p(\Sph)$. Moreover,
$$\|f\|_{\BT^p}=|f(0)|+\|\Q f\|_{L^p(\Sph)}.$$

The following lemma gives some growth estimates for functions in $\BT^p$, which plays a fundamental role in this paper.

\begin{lemma}\label{lem:pointwise}
Let $0<p<\infty$ and $f\in\BT^p$. Then for any $z\in\B$,
\begin{equation}\label{deri-growth}
|Rf(z)|\less\frac{\norm f_{\BT^p}}{(1-|z|^2)^{1+n/p}},
\end{equation}
and
\begin{equation}\label{growth}
|f(z)|\less\frac{\norm f_{\BT^p}}{(1-|z|^2)^{n/p}}
\end{equation}
with implicit constants independent of $f$ and $z$. Moreover, for $0<t\leq1/2$ and $z\in\B$,
\begin{equation}\label{Rftz}
|Rf(tz)|\less t\norm f_{\BT^p}.
\end{equation}
\end{lemma}
\begin{proof}
The estimate \eqref{deri-growth} was established in \cite[Lemma 4.2]{Ch23}, and the equivalence between \eqref{deri-growth} and \eqref{growth} is well-known. Since $Rf(0)=0$, assume $z\neq0$ and define
$$F_z(\lambda):=Rf(\lambda z/2),\quad \lambda\in\D.$$
Then $F_z$ is analytic on $\D$, $F_z(0)=0$, and for any $\lambda\in\D$, the estimate \eqref{deri-growth} gives
$$|F_z(\lambda)|=|Rf(\lambda z/2)|\less\|f\|_{\BT^p}.$$
Hence by the Schwarz lemma, for any $\lambda\in\D$,
$$|Rf(\lambda z/2)|=|F_z(\lambda)|\lesssim|\lambda|\|f\|_{\BT^p}.$$
In particular, for any $0<t\leq1/2$, $|Rf(tz)|\less t\|f\|_{\BT^p}$.
\end{proof}

We also need the following inequality involving the geometry of Kor\'anyi regions $\Gamma(\xi)$.

\begin{lemma}\label{kernel-estimate}
Let $a\in\B$, $\xi\in\Sph$, and let $z\in\Gamma(\xi)$. Then
$$|1-\ip{a}{\xi}|\leq9|1-\ip{z}{a}|.$$
\end{lemma}
\begin{proof}
Write
\[
 \delta=1-|z|^2,\quad d=|1-\ip a\xi|,\quad
 D=|1-\ip za|.
\]
Put $u=\ip z\xi$, $\lambda=\ip a\xi$, and decompose
$$z=u\xi+z_0,\quad a=\lambda\xi+a_0,$$
with $z_0,a_0\perp\xi$. Since $z\in\Gamma(\xi)$, $|1-u|<\delta$. Consequently,
$$|u|\geq1-|1-u|>1-\delta=|z|^2.$$
It follows that
$$|z_0|^2=|z|^2-|u|^2<|z|^2-|z|^4<1-|z|^2=\delta.$$
Moreover,
$$|a_0|^2=|a|^2-|\lambda|^2<1-|\lambda|^2<2|1-\lambda|=2d.$$
Then by the triangle inequality,
\begin{align*}
d&\leq|1-\ip {a}{z}|+|\ip{a}{z}-\ip a\xi|\\
&\leq D+|\lambda||1-u|+|\ip{z_0}{a_0}|\\
&<D+\delta+\sqrt{2d\delta}.
\end{align*}
Noting that $\delta=1-|z|^2<2(1-|z|)\leq2D$, we obtain
\[
 d\leq3D+2\sqrt{Dd}.
\]
Since $D>0$, with $x=\sqrt{d/D}$, this becomes
$x^2\leq3+2x$, whence $x\leq3$ and
$$d\leq9D,$$
as desired.
\end{proof}

The following lemma establishes some necessary conditions for the boundedness of the operators $J_g$, $I_g$ and $M_g$.

\begin{lemma}\label{necessary}
Let $0<p,q<\infty$ and $g\in\hol(\B)$.
\begin{enumerate}
	\item [(1)] If $J_g:\BT^p\to\BT^q$ is bounded, then $g\in\bloch_{1+n/q-n/p}$, and
	$$\|g-g(0)\|_{\bloch_{1+n/q-n/p}}
	\less\|J_g\|_{\BT^p\to\BT^q}.$$
	\item [(2)] If $I_g:\BT^p\to\BT^q$ is bounded, then
	$$\sup_{z\in\B}|g(z)|(1-|z|^2)^{n/q-n/p}
	\less\norm{I_g}_{\BT^p\to\BT^q}.$$
	\item [(3)] If $M_g:\BT^p\to\BT^q$ is bounded, then
	$$\sup_{z\in\B}|g(z)|(1-|z|^2)^{n/q-n/p}
	\less\norm{M_g}_{\BT^p\to\BT^q}.$$
\end{enumerate}
\end{lemma}
\begin{proof}
Fix $\theta>0$ and $a\in\B$. Define
\begin{equation*}
 f_a(z)=\frac{(1-|a|^2)^\theta}
                 {(1-\ip za)^{\theta+n/p}},\quad z\in\B.
\end{equation*}
Then
\[
 Rf_a(z)=
 \frac{(\theta+n/p)(1-|a|^2)^\theta\ip za}
      {(1-\ip za)^{\theta+n/p+1}}.
\]
Using $1-|z|^2\leq2|1-\ip za|$ and Lemma \ref{kernel-estimate}, we get for any $\xi\in\Sph$,
\[
 \Q{f_a}(\xi)=\sup_{z\in\Gamma(\xi)}|Rf_a(z)|(1-|z|^2)\less
 \frac{(1-|a|^2)^\theta}{|1-\ip a\xi|^{\theta+n/p}}.
\]
The standard sphere integral estimate (see, for instance, \cite[Theorem 1.12]{Zh})
\[
 \int_{\Sph}\frac{d\sigma(\xi)}{|1-\ip a\xi|^{n+\theta p}}
 \less(1-|a|^2)^{-\theta p}
\]
therefore gives $f_a\in\BT^p$ with $\norm{f_a}_{\BT^p}\less1$.

(1) If $J_g:\BT^p\to\BT^q$ is bounded, then by \eqref{deri-growth},
\begin{align*}
 \norm{J_g}_{\BT^p\to\BT^q}&\gtrsim\norm{J_gf_a}_{\BT^q}\\
 &\gtrsim|R(J_gf_a)(a)|(1-|a|^2)^{1+n/q}\\
 &=|f_a(a)Rg(a)|(1-|a|^2)^{1+n/q}\\
 &=|Rg(a)|(1-|a|^2)^{1+n/q-n/p}.
\end{align*}
Since $a\in\B$ is arbitrary, we obtain $g\in\bloch_{1+n/q-n/p}$, and
$$\|g-g(0)\|_{\bloch_{1+n/q-n/p}}\less\|J_g\|_{\BT^p\to\BT^q}.$$

(2) If $I_g:\BT^p\to\BT^q$ is bounded, then as above,
\begin{align*}
\norm{I_g}_{\BT^p\to\BT^q}&\gtrsim\norm{I_gf_a}_{\BT^q}\\
&\gtrsim|R(I_gf_a)(a)|(1-|a|^2)^{1+n/q}\\
&=|Rf_a(a)g(a)|(1-|a|^2)^{1+n/q}\\
&\asymp|a|^2|g(a)|(1-|a|^2)^{n/q-n/p}.
\end{align*}
The arbitrariness of $a\in\B$ together with the maximum modulus principle then gives
$$\sup_{z\in\B}|g(z)|(1-|z|^2)^{n/q-n/p}\less \norm{I_g}_{\BT^p\to\BT^q}.$$

(3) If $M_g:\BT^p\to\BT^q$ is bounded, then using \eqref{growth} instead of \eqref{deri-growth} as above, we obtain
\begin{equation*}
\sup_{z\in\B}|g(z)|(1-|z|^2)^{n/q-n/p}\less \norm{M_g}_{\BT^p\to\BT^q}.\qedhere
\end{equation*}
\end{proof}

\section{Proofs of Theorems \ref{main1} and \ref{main2}}\label{proofs of 12}

In this section, we are going to prove Theorems \ref{main1} and \ref{main2}. To this end, we need the boundedness of some fractional maximal functions. Given $\xi,\zeta\in\Sph$, define
$$d(\xi,\zeta):=|1-\ip\xi\zeta|.$$
Then $d(\cdot,\cdot)$ is a quasi-metric on the sphere $\Sph$. For $\xi\in\Sph$ and $\rho>0$, the quasi-metric ball $Q(\xi,\rho)$ is given by
$$Q(\xi,\rho):=\{\zeta\in\Sph:d(\xi,\zeta)<\rho\}.$$
It is well-known that the sphere $\Sph$ with the quasi-distance $d(\cdot,\cdot)$ and the surface measure $d\sigma$ is a space of
homogeneous type, with $\sigma(Q(\xi,\rho))\asymp \rho^n$ for $0<\rho\leq2$.
For $0\leq s<n$, the fractional maximal operator $\M_s$ is defined for measurable functions $F$ on $\Sph$ by
\[
\M_sF(\xi):=\sup_{\rho>0}
\frac{1}{\sigma\big(Q(\xi,\rho)\big)^{1-s/n}}\int_{Q(\xi,\rho)}|F|\dd\sigma,\quad \xi\in\Sph.
\]
We write $\M:=\M_0$. Write $I(z):=\{\xi\in\Sph:z\in\Gamma(\xi)\}$ for $z\in\B\setminus\{0\}$, and write $I(0):=\Sph$. Then $\sigma(I(z))\asymp(1-|z|^2)^n$.

The following boundedness theorem of fractional maximal operators was proved in \cite[Theorem 1.1]{CS}.

\begin{theorem}\label{maximal}
	Let $1<p\leq q<\infty$ and $s=n(1/p-1/q)$. Then $\M_s:L^p(\Sph)\to L^q(\Sph)$ is bounded.
\end{theorem}

We will need the following technical lemma.

\begin{lemma}\label{cap}
Let $\xi\in\Sph$, $z\in\Gamma(\xi)$ and $0\leq t\leq1$. Then
\[
 I(tz)\subset Q\bigl(\xi,8(1-t^2|z|^2)\bigr).
\]
\end{lemma}
\begin{proof}
If $t=0$ or $z=0$, then there is nothing to prove. Assume $t>0$, $z\neq0$, and write $\zeta=z/|z|$. For any $\eta\in I(tz)$, we have
$$|1-t\ip{z}{\eta}|<1-t^2|z|^2,$$
and consequently,
$$|1-\ip{\zeta}{\eta}|\leq|1-t\ip{z}{\eta}|+1-t|z|
<2(1-t^2|z|^2).$$
Since $z\in\Gamma(\xi)$,
$$|1-\ip{\zeta}{\xi}|\leq|1-\ip{z}{\xi}|+1-|z|<2(1-|z|^2)
\leq2(1-t^2|z|^2).$$
The quasi-triangle inequality
$d(\eta,\xi)\leq2\big(d(\eta,\zeta)+d(\zeta,\xi)\big)$
now proves the assertion.
\end{proof}

The following two theorems are the main results of this section, which establish some embedding theorems for the spaces $\BT^p$.

\begin{theorem}\label{embedding}
Let $0<p<q<\infty$ and $s=n(1/p-1/q)$. Then for any $f\in\BT^p$, the function
$$\xi\mapsto\sup_{z\in\Gamma(\xi)}|f(z)|(1-|z|^2)^s$$
belongs to $L^q(\Sph)$, and
\[
 \int_{\Sph}
 \left(\sup_{z\in\Gamma(\xi)}|f(z)|(1-|z|^2)^s\right)^q
 \dd\sigma(\xi)
 \less\norm f_{\BT^p}^q.
\]
\end{theorem}
\begin{proof}
Choose $0<r<p$, and fix $f\in\BT^p$. For any $z\in\B$ and $\eta\in I(z)$, by the definition of $\Q f$,
$$\Q f(\eta)\geq|Rf(z)|(1-|z|^2),$$
which implies
\begin{equation}\label{average}
 |Rf(z)|^r\leq\frac1{(1-|z|^2)^r\sigma(I(z))}
             \int_{I(z)}(\Q f)^r\dd\sigma.
\end{equation}
Fix $\xi\in\Sph$, and let $z\in\Gamma(\xi)$ with $|z|\geq1/2$. Then for any $1/2\leq t\leq1$, Lemma \ref{cap} together with \eqref{average} gives
\begin{align*}
|Rf(tz)|^r&\leq\frac{1}{(1-t^2|z|^2)^r\sigma(I(tz))}
  \int_{I(tz)}(\Q f)^r\dd\sigma\\
&\less\frac{1}{(1-t^2|z|^2)^r\sigma\big(Q(\xi,8(1-t^2|z|^2))\big)}
  \int_{Q(\xi,8(1-t^2|z|^2))}(\Q f)^r\dd\sigma\\
&\lesssim(1-t^2|z|^2)^{-r-sr}\M_{sr}\big((\Q f)^r\big)(\xi),
\end{align*}
which implies
$$|Rf(tz)|\lesssim
(1-t^2|z|^2)^{-1-s}
\bigl[\M_{sr}\big((\Q f)^r\big)(\xi)\bigr]^{1/r}.$$
Integrating over $[1/2,1]$, we obtain
\[
 \int_{1/2}^1|Rf(tz)|\frac{dt}{t}
 \less (1-|z|^2)^{-s}\bigl[\M_{sr}\big((\Q f)^r\big)(\xi)\bigr]^{1/r}
\]
since $s>0$. The estimate \eqref{Rftz} yields
$$\int_{0}^{1/2}|Rf(tz)|\frac{dt}{t}\less\|f\|_{\BT^p},$$
and consequently,
$$|f(z)|\leq\left|\int_0^1Rf(tz)\frac{dt}{t}\right|+|f(0)|
\less (1-|z|^2)^{-s}\bigl[\M_{sr}\big((\Q f)^r\big)(\xi)\bigr]^{1/r}
  +\|f\|_{\BT^p}.$$
Note that the estimate \eqref{growth} implies $\sup_{|z|\leq1/2}|f(z)|\less\|f\|_{\BT^p}$. Therefore,
\[
 \sup_{z\in\Gamma(\xi)}|f(z)|(1-|z|^2)^s
 \less\bigl[\M_{sr}\big((\Q f)^r\big)(\xi)\bigr]^{1/r}+\norm f_{\BT^p}.
\]
Noting that $1<p/r<q/r$ and
$$0<sr=n\left(\frac{1}{p/r}-\frac{1}{q/r}\right)<n,$$
we apply Theorem \ref{maximal} to obtain
\begin{align*}
\left\|\xi\mapsto\sup_{z\in\Gamma(\xi)}
  |f(z)|(1-|z|^2)^s\right\|_{L^q(\Sph)}
&\less\left\|\bigl[\M_{sr}\big((\Q f)^r\big)\bigr]^{1/r}\right\|_{L^q(\Sph)}
  +\|f\|_{\BT^p}\\
&=\left\|\M_{sr}\big((\Q f)^r\big)\right\|_{L^{q/r}(\Sph)}^{1/r}
  +\|f\|_{\BT^p}\\
&\less\|(\Q f)^r\|_{L^{p/r}(\Sph)}^{1/r}+\|f\|_{\BT^p}\\
&=\|\Q f\|_{L^p(\Sph)}+\|f\|_{\BT^p}\\
&\less\|f\|_{\BT^p},
\end{align*}
which establishes the desired estimate.
\end{proof}

\begin{theorem}\label{logembedding}
Let $0<p<\infty$. Then for any $f\in\BT^p$, the function
$$\xi\mapsto\sup_{z\in\Gamma(\xi)}|f(z)|
\left(\log\frac{2}{1-|z|^2}\right)^{-1}$$
belongs to $L^p(\Sph)$, and
	\[
	\int_{\Sph}
	\left(\sup_{z\in\Gamma(\xi)}
	|f(z)|\left(\log\frac{2}{1-|z|^2}\right)^{-1}\right)^p
	\dd\sigma(\xi)\less\norm f_{\BT^p}^p.
	\]
\end{theorem}
\begin{proof}
Choose $0<r<p$. Fix $\xi\in\Sph$ and let $z\in\Gamma(\xi)$.	The proof of Theorem \ref{embedding}, with $s=0$, gives
\[
|Rf(tz)|\less(1-t^2|z|^2)^{-1}\bigl[\M\big((\Q f)^r\big)(\xi)\bigr]^{1/r},
\quad1/2\leq t\leq1.
\]
For $|z|\geq1/2$, integration yields
\[
\int_{1/2}^1|Rf(tz)|\frac{dt}{t}\less
\bigl[\M\big((\Q f)^r\big)(\xi)\bigr]^{1/r}\log\frac 2{1-|z|^2}.
\]
The integral over $(0,1/2)$ and points $|z|<1/2$ are handled by Lemma \ref{lem:pointwise} as in the proof of Theorem \ref{embedding}. Hence we obtain
$$\sup_{z\in\Gamma(\xi)}|f(z)|\left(\log\frac{2}{1-|z|^2}\right)^{-1}
\less \bigl[\M\big((\Q f)^r\big)(\xi)\bigr]^{1/r}+\|f\|_{\BT^p}.$$
Combining this with the $L^{p/r}$ boundedness of $\mathcal M$ provided by Theorem \ref{maximal} finishes the proof.
\end{proof}

We are now ready to prove Theorems \ref{main1} and \ref{main2}.

\begin{proof}[Proof of Theorem \ref{main1}]
The necessity parts all follow from Lemma \ref{necessary}. Suppose now $g\in\bloch_{1+n/q-n/p}$. Then for any $f\in\BT^p$ and $\xi\in\Sph$,
\begin{align*}
\Q(J_gf)(\xi)&=\sup_{z\in\Gamma(\xi)}|f(z)Rg(z)|(1-|z|^2)\\
&\leq\|g-g(0)\|_{\bloch_{1+n/q-n/p}}
  \sup_{z\in\Gamma(\xi)}|f(z)|(1-|z|^2)^{n/p-n/q}.
\end{align*}
This together with Theorem \ref{embedding} gives
\begin{align*}
\int_{\Sph}\Q({J_gf})(\xi)^q\dd\sigma(\xi)
&\leq\|g-g(0)\|_{\bloch_{1+n/q-n/p}}^{q}\int_{\Sph}
  \left(\sup_{z\in\Gamma(\xi)}|f(z)|(1-|z|^2)^{n/p-n/q}\right)^q
  \dd\sigma(\xi)\\
&\less\|g-g(0)\|_{\bloch_{1+n/q-n/p}}^{q}\|f\|_{\BT^p}^q.
\end{align*}
Therefore, $J_g:\BT^p\to\BT^q$ is bounded, and $\|J_g\|_{\BT^p\to\BT^q}\less\|g-g(0)\|_{\bloch_{1+n/q-n/p}}$.
\end{proof}

\begin{proof}[Proof of Theorem \ref{main2}]
The necessity parts all follow from Lemma \ref{necessary} again. We establish the sufficiency parts below.

(1) Suppose $g\in H^{\infty}$. Then for any $f\in\BT^p$ and $\xi\in\Sph$,
$$\Q({I_gf})(\xi)=\sup_{z\in\Gamma(\xi)}|g(z)Rf(z)|(1-|z|^2)
\leq\|g\|_{H^{\infty}}\Q f(\xi),$$
which immediately implies the boundedness of $I_g$ on $\BT^p$ with the estimate $\|I_g\|_{\BT^p\to\BT^p}\leq\|g\|_{H^{\infty}}$.

(3) Suppose $g\in\bloch_{\log}$. Then for any $f\in\BT^p$ and $\xi\in\Sph$,
$$\Q({J_gf})(\xi)=\sup_{z\in\Gamma(\xi)}|f(z)Rg(z)|(1-|z|^2)
\leq\|g-g(0)\|_{\bloch_{\log}}
  \sup_{z\in\Gamma(\xi)}|f(z)|
  \left(\log\frac{2}{1-|z|^2}\right)^{-1},$$
which, in conjunction with Theorem \ref{logembedding}, implies the boundedness of $J_g$ on $\BT^p$ and the estimate $\|J_g\|_{\BT^p\to\BT^p}\less\|g-g(0)\|_{\bloch_{\log}}$.

(5) Suppose $g\in H^{\infty}\cap\bloch_{\log}$. Then by (1) and (3), both $I_g$ and $J_g$ are bounded on $\BT^p$. The boundedness of $M_g$ follows from the identity \eqref{eq:identities}. Moreover, $\|M_g\|_{\BT^p\to\BT^p}\less\|g\|_{H^{\infty}}+\|g\|_{\bloch_{\log}}.$
\end{proof}

\section{Proof of Theorem \ref{main3}}\label{proof of 3}

In this section, we investigate the boundedness of multipliers and integration operators in the case $p>q$. Before proceeding, we establish a synthesis lemma for the spaces $\BT^p$, which plays an essential role in this section.

Let $\beta(\cdot,\cdot)$ be the Bergman metric on $\B$. Given $z\in\B$ and $\delta>0$, we use $D(z,\delta)$ to denote the Bergman metric ball centered at $z$ with radius $\delta$. For $0<\delta\leq1$, let $Z=\{a_k\}\subset\B\setminus\{0\}$ be a $\delta$-lattice in the Bergman metric of $\B$, that is, a sequence with the following conditions (see \cite[Theorem 2.23]{Zh}):
\begin{enumerate}
	\item [(i)] $\B=\bigcup_kD(a_k,\delta)$;
	\item [(ii)] the sets $D(a_k,\delta/4)$ are pairwise disjoint;
	\item [(iii)] there exists a positive integer $N$, independent of $\delta$, such that each point $z\in\B$ belongs to at most $N$ of the sets $D(a_k,4\delta)$.
\end{enumerate}
Given $0<p<\infty$, the discrete tent space $T_\infty^p(Z)$ consists of
sequences $\lambda=\{\lambda_k\}$ such that
\[
 \norm \lambda_{T_\infty^p(Z)}:=
 \left(\int_{\Sph}
 \left(\sup_{a_k\in\Gamma(\xi)}|\lambda_k|\right)^p
 \dd\sigma(\xi)\right)^{1/p}<\infty.
\]
The supremum over an empty index set is understood to be zero.

\begin{lemma}\label{synthesis}
Let $0<p<\infty$, $\theta>n+n/p$, and let $Z=\{a_k\}\subset\B\setminus\{0\}$ be a $\delta$-lattice. Then for any $\lambda=\{\lambda_k\}\in T^p_{\infty}(Z)$, the function
\begin{equation}\label{systhesis-series}
S(\lambda)(z):=\sum_{k=1}^{\infty}\lambda_k
\left(\frac{1-|a_k|^2}{1-\ip{z}{a_k}}\right)^{\theta}
\end{equation}
belongs to $\BT^p$, and $\|S(\lambda)\|_{\BT^p}\lesssim\|\lambda\|_{T^p_{\infty}(Z)}$ with implicit constant depending on $n,p,\theta$ and $\delta$.
\end{lemma}
\begin{proof}
Since $\theta>n+n/p$, we may choose $0<s<p$ such that $\theta-n/s>n$. For $\lambda=\{\lambda_k\}\in T^p_{\infty}(Z)$, define
$$\Lambda(\xi):=\sup_{a_k\in\Gamma(\xi)}|\lambda_k|,\quad
\xi\in\Sph,$$
and
$$G(\xi):=\bigl[\M(\Lambda^s)(\xi)\bigr]^{1/s},\quad \xi\in\Sph.$$
Since $p/s>1$, Theorem \ref{maximal} gives
\begin{equation}\label{maximal-bound}
 \|G\|_{L^p(\Sph)}^s
 =\|\M(\Lambda^s)\|_{L^{p/s}(\Sph)}
 \lesssim\|\Lambda^s\|_{L^{p/s}(\Sph)}
 =\|\lambda\|_{T_\infty^p(Z)}^s.
\end{equation}

We now show that the series \eqref{systhesis-series} defines a holomorphic function on $\B$. The definition of $\Lambda$ gives that for any $k\geq1$,
$$\|\lambda\|_{T^p_{\infty}(Z)}^p=\|\Lambda\|_{L^p(\Sph)}^p
\geq\int_{I(a_k)}\Lambda(\xi)^p\dd\sigma(\xi)
\gtrsim|\lambda_k|^p(1-|a_k|^2)^n,$$
which implies that for $|z|\leq r<1$,
\[
\left|\lambda_k\left(\frac{1-|a_k|^2}{(1-\ip z{a_k})}\right)^{\theta}\right|
\lesssim (1-r)^{-\theta}\|\lambda\|_{T_\infty^p(Z)}
(1-|a_k|^2)^{\theta-n/p}.
\]
Since $\theta-n/p>n$, we have
\begin{align*}
\sum_{k=1}^{\infty}(1-|a_k|^2)^{\theta-n/p}
&\asymp\sum_{k=1}^{\infty}
\int_{D(a_k,\delta)}(1-|w|^2)^{\theta-n/p-n-1}\dd V(w)\\
&\lesssim\int_{\B}(1-|w|^2)^{\theta-n/p-n-1}\dd V(w)<\infty.
\end{align*}
Therefore, the series \eqref{systhesis-series} converges absolutely and uniformly on compact subsets of $\B$, and so defines a holomorphic function on $\B$. In particular, the same estimate also gives
\[
|S(\lambda)(0)|\leq\sum_{k=1}^{\infty}|\lambda_k|(1-|a_k|^2)^\theta
\lesssim\|\lambda\|_{T_\infty^p(Z)}.
\]

Fix $\xi\in\Sph$ and $z\in\Gamma(\xi)$. We claim that for every $k\geq1$,
\begin{equation}\label{coefficient-bound}
 |\lambda_k|\lesssim
 \left(\frac{|1-\ip z{a_k}|}{1-|a_k|^2}\right)^{n/s}G(\xi).
\end{equation}
Indeed, by the definition of $\Lambda$, $\Lambda\geq|\lambda_k|$ on $I(a_k)$, so
\begin{equation}\label{lambda-mean}
|\lambda_k|^s\leq
\frac1{\sigma(I(a_k))}\int_{I(a_k)}\Lambda^s\dd\sigma.
\end{equation}
We enclose $I(a_k)$ in a cap centered at $\xi$. If $|a_k|\geq1/2$,
put $\zeta_k=a_k/|a_k|$. For $\eta\in I(a_k)$,
\[
 |1-\ip{\zeta_k}\eta|
 \leq |1-\ip{a_k}\eta|+(1-|a_k|)\leq2(1-|a_k|^2)
 \leq4|1-\ip{z}{a_k}|.
\]
Lemma \ref{kernel-estimate} gives
\[
 |1-\ip{\zeta_k}\xi|
 \leq |1-\ip{a_k}\xi|+(1-|a_k|)\leq10|1-\ip{z}{a_k}|.
\]
Consequently,
$$|1-\ip{\eta}{\xi}|
\leq2|1-\ip{\zeta_k}{\eta}|+2|1-\ip{\zeta_k}{\xi}|
\leq28|1-\ip{z}{a_k}|.$$
That is, $I(a_k)\subset Q(\xi,28|1-\ip{z}{a_k}|)$. If $|a_k|<1/2$, then $|1-\ip{z}{a_k}|\geq1/2$,
and we instead use the whole sphere $\Sph$. In either case there is a cap $Q_k$
centered at $\xi$ with
\[
 I(a_k)\subset Q_k,\qquad \sigma(Q_k)\lesssim |1-\ip{z}{a_k}|^n.
\]
Consequently, the inequality \eqref{lambda-mean} yields
\[
 |\lambda_k|^s\lesssim
 \frac{1}{(1-|a_k|^2)^n}\int_{Q_k}\Lambda^s\,d\sigma
 \lesssim\left(\frac{|1-\ip{z}{a_k}|}{1-|a_k|^2}\right)^n
 \M(\Lambda^s)(\xi),
\]
which proves \eqref{coefficient-bound}.

For every positive integer $K$, we next consider the finite sum
$$S_{K}(\lambda)(z):=\sum_{k=1}^{K}\lambda_k
\left(\frac{1-|a_k|^2}{1-\ip z{a_k}}\right)^{\theta}.$$
Differentiation gives
\[
\big|R\big(S_K(\lambda)\big)(z)\big|\lesssim
\sum_{k=1}^K|\lambda_k|\frac{(1-|a_k|^2)^\theta}
  {|1-\ip{z}{a_k}|^{\theta+1}}.
\]
For $z\in\Gamma(\xi)$, \eqref{coefficient-bound} implies
\begin{align*}
\big|R\big(S_K(\lambda)\big)(z)\big|(1-|z|^2)
\lesssim G(\xi)(1-|z|^2)\sum_{k=1}^{\infty}
\frac{(1-|a_k|^2)^{\theta-n/s}}{|1-\ip{z}{a_k}|^{\theta+1-n/s}}.
\end{align*}
Since $\theta-n/s>n$, the classical Forelli--Rudin estimate (see \cite[Theorem 1.12]{Zh}) together with the subharmonic property gives
$$\sum_{k=1}^{\infty}
\frac{(1-|a_k|^2)^{\theta-n/s}}{|1-\ip{z}{a_k}|^{\theta+1-n/s}}
\less\int_{\B}\frac{(1-|w|^2)^{\theta-n/s-n-1}}{|1-\ip{z}{w}|^{\theta+1-n/s}}\dd V(w)
\lesssim\frac{1}{1-|z|^2}.$$
Consequently, for any $\xi\in\Sph$,
$$\Q\big({S_K(\lambda)}\big)(\xi)
=\sup_{z\in\Gamma(\xi)}\big|R\big(S_K(\lambda)\big)(z)\big|(1-|z|^2)
\less G(\xi)$$
with implicit constant independent of $K$. Since $R\big(S_K(\lambda)\big)\to R\big(S(\lambda)\big)$ uniformly on compact subsets of $\B$ as $K\to\infty$, we obtain
$$\Q\big({S(\lambda)}\big)(\xi)
\leq\liminf_{K\to\infty}\Q\big({S_K(\lambda)}\big)(\xi)
\less G(\xi).$$
This together with \eqref{maximal-bound} yields
$$\left\|\Q\big({S(\lambda)}\big)\right\|_{L^p(\Sph)}
\less\|G\|_{L^p(\Sph)}\less\|\lambda\|_{T^p_{\infty}(Z)}.$$
Since $|S(\lambda)(0)|\less\|\lambda\|_{T^p_{\infty}(Z)}$, we conclude that $S(\lambda)\in\BT^p$, and
$$\|S(\lambda)\|_{\BT^p}
=|S(\lambda)(0)|+\left\|\Q\big({S(\lambda)}\big)\right\|_{L^p(\Sph)}
\less\|\lambda\|_{T^p_{\infty}(Z)}.$$
The proof is complete.
\end{proof}

We also need the characterization of bounded multipliers between discrete tent spaces $T^p_{\infty}(Z)$. Given a sequence $b=\{b_k\}$, the pointwise multiplier $\mathfrak{D}_b$ is defined by
$$\mathfrak{D}_b\lambda:=\{b_k\lambda_k\},\quad \lambda=\{\lambda_k\}.$$

\begin{lemma}\label{discrete-tent-multipliers}
Let $0<q<p<\infty$, $Z=\{a_k\}\subset\B\setminus\{0\}$ be a $\delta$-lattice and $b=\{b_k\}$. Then $\mathfrak{D}_b:T_\infty^p(Z)\to T_\infty^q(Z)$ is bounded if and only if $b\in T_\infty^{\frac{pq}{p-q}}(Z)$. Moreover,
\[
 \|\mathfrak{D}_b\|_{T_\infty^p(Z)\to T_\infty^q(Z)}
 =\|b\|_{T_\infty^{pq/(p-q)}(Z)}.
\]
\end{lemma}
\begin{proof}
The sufficiency together with the norm inequality
$$\|\mathfrak{D}_{b}\|_{T^p_{\infty}(Z)\to T^q_{\infty}(Z)}
\leq\|b\|_{T^{pq/(p-q)}_{\infty}(Z)}$$
follows directly from H\"older's inequality.  Suppose now that $\mathfrak{D}_b:T_\infty^p(Z)\to T_\infty^q(Z)$ is bounded. For each positive integer $m$, define $\lambda^{(m)}=\{\lambda^{(m)}_k\}$ by
\begin{equation*}
\lambda_k^{(m)}=
\begin{cases}
	|b_k|^{q/(p-q)},&1\leq k\leq m,\\
	0,&k>m,
\end{cases}
\end{equation*}
and put
\[
 B_m(\xi):=\sup_{\substack{1\leq k\leq m\\ a_k\in\Gamma(\xi)}}|b_k|,\quad \xi\in\Sph.
\]
The sequence $\lambda^{(m)}$ has finite support and belongs to
$T_\infty^p(Z)$. Moreover,
\[
 \sup_{a_k\in\Gamma(\xi)}|\lambda_k^{(m)}|=B_m(\xi)^{q/(p-q)},
\]
and
\[
 \sup_{a_k\in\Gamma(\xi)}|b_k\lambda_k^{(m)}|=B_m(\xi)^{p/(p-q)}.
\]
Applying the assumed boundedness of $\mathfrak{D}_b$ to $\lambda^{(m)}$, we obtain
\[
\left(\int_{\Sph}B_m^{\frac{pq}{p-q}}\dd\sigma\right)^{1/q}
\leq\|\mathfrak{D}_{b}\|_{T^p_{\infty}(Z)\to T^q_{\infty}(Z)}
  \left(\int_{\Sph}B_m^{\frac{pq}{p-q}}\dd\sigma\right)^{1/p}.
\]
The integral is finite because $B_m$ is bounded. Hence we obtain
\[
\left(\int_{\Sph}B_m^{\frac{pq}{p-q}}
  \dd\sigma\right)^{\frac{p-q}{pq}}
\leq\|\mathfrak{D}_{b}\|_{T^p_{\infty}(Z)\to T^q_{\infty}(Z)}.
\]
Letting $m\to\infty$ and applying the monotone convergence theorem, we establish
$$\left(\int_{\Sph}\left(\sup_{a_k\in\Gamma(\xi)}
|b_k|\right)^{\frac{pq}{p-q}}\dd\sigma(\xi)\right)^{\frac{p-q}{pq}}
\leq\|\mathfrak{D}_{b}\|_{T^p_{\infty}(Z)\to T^q_{\infty}(Z)}.$$
This proves the necessity and the inequality $\|b\|_{T^{pq/(p-q)}_{\infty}(Z)}
\leq\|\mathfrak{D}_{b}\|_{T^p_{\infty}(Z)\to T^q_{\infty}(Z)}$.
\end{proof}

The following discretization lemma can be found in \cite[Lemma 3]{MPPW}.

\begin{lemma}\label{discretization}
Let $0<p<\infty$ and $\alpha\geq0$. There exists $\delta_0\in(0,1)$ so that, if $0<\delta<\delta_0$ and $Z=\{a_k\}\subset\B\setminus\{0\}$ is a $\delta$-lattice, then for any holomorphic function $f$ on $\B$ such that $f(\cdot)(1-|\cdot|^2)^{\alpha}\in T^p_{\infty}$,
$$\int_{\Sph}\left(\sup_{z\in\Gamma(\xi)}
  |f(z)|(1-|z|^2)^{\alpha}\right)^p\dd\sigma(\xi)
\less\int_{\Sph}\left(\sup_{a_k\in\Gamma(\xi)}
  |f(a_k)|(1-|a_k|^2)^{\alpha}\right)^p\dd\sigma(\xi)$$
with implicit constant depending only on $p,\alpha,n$.
\end{lemma}

We now consider the multipliers and integration operators induced by dilated symbols. For $0<\rho<1$ and $g\in\hol(\B)$, write $g_{\rho}(z):=g(\rho z)$.

\begin{lemma}\label{dilated-symbol-bound}
Let $0<p,q<\infty$ and $g\in\hol(\B)$. Then:
\begin{enumerate}
	\item [(1)] If $I_g:\BT^p\to\BT^q$ is bounded, then for any $0<\rho<1$, $I_{g_{\rho}}:\BT^p\to\BT^q$ is bounded, and
	\[
	\|I_{g_\rho}\|_{\BT^p\to\BT^q}\leq\|I_g\|_{\BT^p\to\BT^q}.
	\]
	\item [(2)] If $J_g:\BT^p\to\BT^q$ is bounded, then for any $0<\rho<1$, $J_{g_{\rho}}:\BT^p\to\BT^q$ is bounded, and
	\[
	\|J_{g_\rho}\|_{\BT^p\to\BT^q}\leq\|J_g\|_{\BT^p\to\BT^q}.
	\]
	\item [(3)] If $M_g:\BT^p\to\BT^q$ is bounded, then for any $0<\rho<1$, $M_{g_{\rho}}:\BT^p\to\BT^q$ is bounded, and
	\[
	\|M_{g_\rho}\|_{\BT^p\to\BT^q}\leq2\|M_g\|_{\BT^p\to\BT^q}.
	\]
\end{enumerate}
\end{lemma}
\begin{proof}
For $\vartheta\in\mathbb R$ and $f\in\hol(\B)$, define
\[
 U_\vartheta f(z)=f(e^{i\vartheta}z),\quad z\in\B.
\]
The identities
\[
 R(U_\vartheta f)(z)=Rf(e^{i\vartheta}z)\qquad\text{and}\qquad
 e^{i\vartheta}\Gamma(\xi)=\Gamma(e^{i\vartheta}\xi)
\]
imply $\Q({U_\vartheta f})(\xi)=\Q f(e^{i\vartheta}\xi)$ for any $\xi\in\Sph$. Since the surface measure on $\Sph$ is rotation invariant and $U_\vartheta f(0)=f(0)$, $U_\vartheta$ is an isometry on $\BT^t$ for every $0<t<\infty$.

(1) The definition of $I_g$ easily gives
\[
 I_{U_{\vartheta}g}=U_\vartheta I_gU_{-\vartheta}.
\]
Consequently, for any $\vartheta\in\mathbb{R}$,
\begin{equation}\label{rotated-symbol-norm}
 \|I_{U_{\vartheta}g}\|_{\BT^p\to\BT^q}=\|I_g\|_{\BT^p\to\BT^q}.
\end{equation}
Fix $0<s<q$ and $0<\rho<1$. Let
\[
 P_\rho(\vartheta):=\frac{1-\rho^2}{|e^{i\vartheta}-\rho|^2}
\]
be the Poisson kernel of the unit disk $\D$. For each fixed $z\in\B$, the function $\lambda\mapsto g(\lambda z)$ is analytic on a neighborhood of $\overline{\D}$. Then by the subharmonic property,
\[
 |g(\rho z)|^s\leq\int_0^{2\pi}P_\rho(\vartheta)
 |g(e^{i\vartheta}z)|^s\,\frac{d\vartheta}{2\pi}.
\]
For $f\in\BT^p$, multiplying the above inequality by
$\big(|Rf(z)(1-|z|^2)|\big)^s$ and applying the identities $R(I_{g_\rho}f)=g_\rho Rf$ and $R(I_{U_\vartheta g}f)=(U_\vartheta g) Rf$, we obtain
$$\big(|R(I_{g_\rho}f)(z)|(1-|z|^2)\big)^s\leq
\int_0^{2\pi}P_\rho(\vartheta)
\big(|R(I_{U_{\vartheta}g}f)(z)|(1-|z|^2)\big)^s\,\frac{d\vartheta}{2\pi},$$
which implies that for any $\xi\in\Sph$,
\begin{equation*}
	\Q({I_{g_\rho}f})(\xi)^s
	\leq\int_0^{2\pi}P_\rho(\vartheta)
	\Q({I_{U_\vartheta g}f})(\xi)^s\,\frac{d\vartheta}{2\pi}.
\end{equation*}
Since $q/s>1$, Minkowski's inequality together with \eqref{rotated-symbol-norm} gives
\begin{align*}
 \|\Q({I_{g_\rho}f})\|_{L^q(\Sph)}^s
 &=\|\Q({I_{g_\rho}f})^s\|_{L^{q/s}(\Sph)}\\
 &\leq\int_0^{2\pi}P_\rho(\vartheta)
 \|\Q({I_{U_\vartheta g}f})^s\|_{L^{q/s}(\Sph)}\,\frac{d\vartheta}{2\pi}\\
 &=\int_0^{2\pi}P_\rho(\vartheta)
 \|\Q({I_{U_\vartheta g}f})\|_{L^q(\Sph)}^s\,\frac{d\vartheta}{2\pi}\\
 &\leq\|I_g\|_{\BT^p\to\BT^q}^s\|f\|_{\BT^p}^s.
\end{align*}
Since $I_{g_\rho}f(0)=0$, we conclude that
\[
 \|I_{g_\rho}f\|_{\BT^q}\leq \|I_g\|_{\BT^p\to\BT^q}\|f\|_{\BT^p}.
\]

(2) This is exactly the same as (1). In fact, $J_{U_{\vartheta}g}=U_{\vartheta}J_gU_{-\vartheta}$. We can use the subharmonic property of the function $\lambda\mapsto|Rg(\lambda z)|^s$ and argue as in (1) to establish the desired result.

(3) As in (1), $M_{U_{\vartheta}g}=U_{\vartheta}M_gU_{-\vartheta}$. For fixed $f\in\BT^p$ and $z\in\B$, using the subharmonic property of the function
$$\lambda\mapsto|g(\lambda z)Rf(z)+f(z)Rg(\lambda z)|^s,$$
we obtain
$$|R(g_{\rho}f)(z)|^s\leq\int_0^{2\pi}P_{\rho}(\vartheta)
|R(M_{U_{\vartheta}g}f)(z)|^s\,\frac{d\vartheta}{2\pi}.$$
As in (1), this implies
$$\|\Q({g_{\rho}f})\|_{L^q(\Sph)}\leq
\|M_g\|_{\BT^p\to\BT^q}\|f\|_{\BT^p}.$$
Testing $M_g$ on the constant function $1$ easily gives $|g(0)|\leq\|M_g\|_{\BT^p\to\BT^q}$. Therefore,
\begin{equation*}
\|M_{g_{\rho}}f\|_{\BT^q}
=|g(0)f(0)|+\|\Q({g_{\rho}f})\|_{L^q(\Sph)}
\leq2\|M_g\|_{\BT^p\to\BT^q}\|f\|_{\BT^p}.\qedhere
\end{equation*}
\end{proof}

The following theorem characterizes the boundedness of $I_g:\BT^p\to\BT^q$ in the case $p>q$.

\begin{theorem}\label{I_gp>q}
Let $0<q<p<\infty$ and $g\in\hol(\B)$. Then $I_g:\BT^p\to\BT^q$ is bounded if and only if $g\in H^{\frac{pq}{p-q}}$. Moreover,
$$\|I_g\|_{\BT^p\to\BT^q}\asymp\|g\|_{H^{pq/(p-q)}}.$$
\end{theorem}
\begin{proof}
Suppose first that $g\in H^{\frac{pq}{p-q}}$. Recall that the non-tangential maximal function $\N g$ is defined by
$$\N g(\xi):=\sup_{z\in\Gamma(\xi)}|g(z)|,\quad \xi\in\Sph.$$
Then for any $f\in\BT^p$, $\Q({I_gf})\leq \Q f\cdot \N g$, and H\"older's inequality gives
\[
\norm{\Q({I_gf})}_{L^q(\Sph)}\leq\|\Q f\cdot\N g\|_{L^q(\Sph)}
\leq\norm{\Q f}_{L^p(\Sph)}\norm{\N g}_{L^{\frac{pq}{p-q}}(\Sph)}
\less\norm f_{\BT^p}\norm g_{H^{\frac{pq}{p-q}}}.
\]
The final step is the standard non-tangential maximal estimate for Hardy spaces; see \cite[Theorem 4.24]{Zh}. Therefore, $I_g:\BT^p\to\BT^q$ is bounded, and $\|I_g\|_{\BT^p\to\BT^q}\lesssim\|g\|_{H^{\frac{pq}{p-q}}}$.

Conversely, suppose that $I_g:\BT^p\to\BT^q$ is bounded. In view of Lemma \ref{dilated-symbol-bound}, it suffices to establish the norm estimate 
$$\|g\|_{H^{pq/(p-q)}}\less\|I_g\|_{\BT^p\to\BT^q}$$
in the case $g\in H^{\frac{pq}{p-q}}$. If $g\in H^{\frac{pq}{p-q}}$, then the non-tangential maximal estimate for Hardy spaces implies $g\in T^{\frac{pq}{p-q}}_{\infty}$. By Lemma \ref{discretization}, it is sufficient to prove that
\begin{equation}\label{aim-dis}
\int_{\Sph}\left(\sup_{a_k\in\Gamma(\xi)}|g(a_k)|\right)
^{\frac{pq}{p-q}}\dd\sigma(\xi)
\less\|I_g\|_{\BT^p\to\BT^q}^{\frac{pq}{p-q}}
\end{equation}
when $Z=\{a_k\}\subset\B\setminus\{0\}$ is a $\delta$-lattice and $\delta$ is small enough.

Let $Z=\{a_k\}\subset\B\setminus\{0\}$ be a $\delta$-lattice with sufficiently small $\delta>0$. For any $\lambda=\{\lambda_k\}\in T^p_{\infty}(Z)$, consider the test functions
\begin{equation}\label{testp>q}
S_{\tau}(\lambda)(z)
:=\sum_{k=1}^{\infty}\lambda_kr_k(\tau)f_k(z),
\end{equation}
where $\{r_k\}$ is a sequence of Rademacher functions on $[0,1]$, and
$$f_k(z)=\left(\frac{1-|a_k|^2}{1-\ip{z}{a_k}}\right)^{\theta}$$
with $\theta>n+n/p$. By Lemma \ref{synthesis}, the random series \eqref{testp>q} and its derivatives converge absolutely and uniformly on compact subsets of $\B$, uniformly in the choices of signs. Moreover, $S_{\tau}(\lambda)\in\BT^p$ for almost every $\tau\in[0,1]$, and $\|S_{\tau}(\lambda)\|_{\BT^p}\less\|\lambda\|_{T^p_{\infty}(Z)}$. Then by the boundedness of $I_g:\BT^p\to\BT^q$, we have for almost every $\tau\in[0,1]$,
$$\|I_g\|_{\BT^p\to\BT^q}^q\|\lambda\|_{T^p_{\infty}(Z)}^q
\gtrsim\int_{\Sph}\left(\sup_{z\in\Gamma(\xi)}
  \left|\sum_{k=1}^{\infty}\lambda_kr_k(\tau)Rf_k(z)\right|
  |g(z)|(1-|z|^2)\right)^q\dd\sigma(\xi).$$
Integrating with respect to $\tau$ and using Fubini's theorem and Khinchine's inequality, we get
\begin{align*}
&\|I_g\|_{\BT^p\to\BT^q}^q\|\lambda\|_{T^p_{\infty}(Z)}^q\\
&\ \gtrsim\int_{\Sph}\int_0^1\left(\sup_{z\in\Gamma(\xi)}
  \left|\sum_{k=1}^{\infty}\lambda_kr_k(\tau)Rf_k(z)
  g(z)(1-|z|^2)\right|^q\right)\dd\tau\dd\sigma(\xi)\\
&\ \geq\int_{\Sph}\left(\sup_{z\in\Gamma(\xi)}\int_0^1
  \left|\sum_{k=1}^{\infty}\lambda_kr_k(\tau)Rf_k(z)
  g(z)(1-|z|^2)\right|^q\dd\tau\right)\dd\sigma(\xi)\\
&\ \asymp\int_{\Sph}\left(\sup_{z\in\Gamma(\xi)}
  \left(\sum_{k=1}^{\infty}|\lambda_k|^2|Rf_k(z)|^2
  |g(z)|^2(1-|z|^2)^2\right)^{q/2}
  \right)\dd\sigma(\xi)\\
&\ \geq\int_{\Sph}\left(\sup_{a_j\in\Gamma(\xi)}
  \left(\sum_{k=1}^{\infty}|\lambda_k|^2|Rf_k(a_j)|^2
  |g(a_j)|^2(1-|a_j|^2)^2\right)^{q/2}
  \right)\dd\sigma(\xi)\\
&\ \geq\int_{\Sph}\left(\sup_{a_j\in\Gamma(\xi)}
  |\lambda_j||Rf_j(a_j)||g(a_j)|(1-|a_j|^2)
  \right)^{q}\dd\sigma(\xi).
\end{align*}
Noting that
$$Rf_j(a_j)=\theta\frac{|a_j|^2}{1-|a_j|^2},$$
we obtain
$$\|I_g\|_{\BT^p\to\BT^q}^q\|\lambda\|_{T^p_{\infty}(Z)}^q
\gtrsim\int_{\Sph}\left(\sup_{a_j\in\Gamma(\xi)}
  |\lambda_j||a_j|^2|g(a_j)|\right)^q\dd\sigma(\xi).$$
That is, the sequence $\{\lambda_k|a_k|^2g(a_k)\}$ belongs to $T^q_{\infty}(Z)$, and
$$\left\|\{\lambda_k|a_k|^2g(a_k)\}\right\|_{T^q_{\infty}(Z)}
\less\|I_g\|_{\BT^p\to\BT^q}\|\lambda\|_{T^p_{\infty}(Z)}.$$
Since $\lambda\in T^p_{\infty}(Z)$ is arbitrary, we apply Lemma \ref{discrete-tent-multipliers} to establish
\begin{equation}\label{aim-dis'}
\int_{\Sph}\left(\sup_{a_k\in\Gamma(\xi)}|a_k|^2|g(a_k)|\right)
^{\frac{pq}{p-q}}\dd\sigma(\xi)
\less\|I_g\|_{\BT^p\to\BT^q}^{\frac{pq}{p-q}}.
\end{equation}
To deduce \eqref{aim-dis} from \eqref{aim-dis'}, we separate the
lattice points near the origin from those near the boundary.
By Lemma \ref{necessary},
\[
\sup_{z\in\B}
|g(z)|(1-|z|^2)^{n/q-n/p}
\lesssim \|I_g\|_{\BT^p\to\BT^q}.
\]
In particular,
\[
\sup_{|z|\leq1/2}|g(z)|
\lesssim \|I_g\|_{\BT^p\to\BT^q}.
\]
On the other hand, if $|a_k|>1/2$, then
$|g(a_k)|\leq4|a_k|^2|g(a_k)|$. Consequently,
\begin{align*}
&\int_{\Sph}
  \left(\sup_{a_k\in\Gamma(\xi)}|g(a_k)|\right)^{\frac{pq}{p-q}}
  \dd\sigma(\xi)\\
&\ \lesssim\|I_g\|_{\BT^p\to\BT^q}^{\frac{pq}{p-q}}+
	\int_{\Sph}\left(\sup_{a_k\in\Gamma(\xi)}
	|a_k|^2|g(a_k)|\right)^{\frac{pq}{p-q}}
	\dd\sigma(\xi)\\
&\ \lesssim \|I_g\|_{\BT^p\to\BT^q}^{\frac{pq}{p-q}},
\end{align*}
where the last inequality follows from \eqref{aim-dis'}. This proves \eqref{aim-dis} and finishes the proof.
\end{proof}

The following theorem gives separate necessary and sufficient conditions for the boundedness of $J_g:\BT^p\to\BT^q$ in the case $p>q$.

\begin{theorem}\label{J_gp>q}
Let $0<q<p<\infty$ and $g\in\hol(\B)$. Then:
\begin{enumerate}
	\item [(1)] If $J_g:\BT^p\to\BT^q$ is bounded, then $g\in\BT^{\frac{pq}{p-q}}$, and
	$$\|g-g(0)\|_{\BT^{pq/(p-q)}}\less\|J_g\|_{\BT^p\to\BT^q}.$$
	\item [(2)] If $g\in\BT^{\frac{pq}{p-q}}_{\log}$, then $J_g:\BT^p\to\BT^q$ is bounded, and
	$$\|J_g\|_{\BT^p\to\BT^q}
	\less\|g-g(0)\|_{\BT^{pq/(p-q)}_{\log}}.$$
\end{enumerate}
\end{theorem}
\begin{proof}
(1) This is similar to the necessity part of Theorem \ref{I_gp>q}. In fact, Lemma \ref{dilated-symbol-bound} ensures that it suffices to establish the norm estimate
$$\|g-g(0)\|_{\BT^{pq/(p-q)}}\less\|J_g\|_{\BT^p\to\BT^q}$$
in the case $g\in\BT^{\frac{pq}{p-q}}$. By Lemma \ref{discretization}, it is sufficient to prove
$$\int_{\Sph}\left(\sup_{a_k\in\Gamma(\xi)}|Rg(a_k)|(1-|a_k|^2)\right)
^{\frac{pq}{p-q}}\dd\sigma
\less\|J_g\|_{\BT^p\to\BT^q}^{\frac{pq}{p-q}}$$
when $Z=\{a_k\}\subset\B\setminus\{0\}$ is a $\delta$-lattice and $\delta$ is small enough. Using the test functions \eqref{testp>q}, we have
$$\|J_g\|_{\BT^p\to\BT^q}^q\|\lambda\|_{T^p_{\infty}(Z)}^q
\gtrsim\int_{\Sph}\sup_{z\in\Gamma(\xi)}
\left|\sum_{k=1}^{\infty}\lambda_kr_k(\tau)f_k(z)Rg(z)(1-|z|^2)\right|^q\dd\sigma(\xi).$$
Arguing as in the proof of Theorem \ref{I_gp>q}, this implies
$$\|J_g\|_{\BT^p\to\BT^q}^q\|\lambda\|_{T^p_{\infty}(Z)}^q
\gtrsim\int_{\Sph}\left(\sup_{a_j\in\Gamma(\xi)}|\lambda_j|
  |Rg(a_j)|(1-|a_j|^2)\right)^q\dd\sigma(\xi),$$
which, together with Lemma \ref{discrete-tent-multipliers}, gives the desired estimate.

(2) If $g\in\BT^{\frac{pq}{p-q}}_{\log}$, then for any $f\in\BT^p$, H\"older's inequality yields
\begin{align*}
&\int_{\Sph}\left(\sup_{z\in\Gamma(\xi)}|R(J_gf)(z)|
  (1-|z|^2)\right)^q\dd\sigma(\xi)\\
&\ \leq\int_{\Sph}\left(\sup_{z\in\Gamma(\xi)}
  |Rg(z)|(1-|z|^2)\log\frac{2}{1-|z|^2}\right)^q\cdot
  \left(\sup_{z\in\Gamma(\xi)}\frac{|f(z)|}{\log\frac{2}{1-|z|^2}}
  \right)^q\dd\sigma(\xi)\\
&\ \leq\|g-g(0)\|_{\BT^{pq/(p-q)}_{\log}}^q
  \left(\int_{\Sph}\left(\sup_{z\in\Gamma(\xi)}
  \frac{|f(z)|}{\log\frac{2}{1-|z|^2}}\right)^p\dd\sigma(\xi)
  \right)^{q/p}\\
&\ \less\|g-g(0)\|_{\BT^{pq/(p-q)}_{\log}}^q\|f\|_{\BT^p}^q,
\end{align*}
where the last step is due to Theorem \ref{logembedding}. Since $J_gf(0)=0$, we conclude that $J_g:\BT^p\to\BT^q$ is bounded, and
\begin{equation*}
\|J_g\|_{\BT^p\to\BT^q}
\less\|g-g(0)\|_{\BT^{pq/(p-q)}_{\log}}.\qedhere
\end{equation*}
\end{proof}

The following theorem considers the boundedness of $M_g:\BT^p\to\BT^q$ in the case $p>q$.

\begin{theorem}\label{M_gp>q}
Let $0<q<p<\infty$ and $g\in\hol(\B)$. Then
\begin{enumerate}
\item [(1)] If $g\in 
            H^{\frac{pq}{p-q}}\cap\BT_{\log}^{\frac{pq}{p-q}}$, then $M_g:\BT^p\to\BT^q$ is bounded, and
            $$\|M_g\|_{\BT^p\to\BT^q}
            \lesssim\|g\|_{H^{pq/(p-q)}}
            +\|g\|_{\BT_{\log}^{pq/(p-q)}}.$$
\item [(2)] If $M_g:\BT^p\to\BT^q$ is bounded, then
            $g\in H^{\frac{pq}{p-q}}$, and
            $$\|g\|_{H^{pq/(p-q)}}\lesssim
            \|M_g\|_{\BT^p\to\BT^q}.$$
\item [(3)] $M_g:\BT^p\to\BT^q$ is bounded if and only if $g\in 
            H^{\frac{pq}{p-q}}$ and $J_g:\BT^p\to\BT^q$ is bounded. In this case,
            $$\|M_g\|_{\BT^p\to\BT^q}
            \asymp\|g\|_{H^{pq/(p-q)}}+\|J_g\|_{\BT^p\to\BT^q}.$$
\end{enumerate}
\end{theorem}
\begin{proof}
(1) If $g\in H^{\frac{pq}{p-q}}\cap\BT_{\log}^{\frac{pq}{p-q}}$, then by Theorems \ref{I_gp>q} and \ref{J_gp>q}, both $I_g$ and $J_g$ are bounded from $\BT^p$ into $\BT^q$. Hence the boundedness of $M_g:\BT^p\to\BT^q$ together with the desired norm estimate follows from the identity \eqref{eq:identities}.

(2) As in the proof of Theorem \ref{I_gp>q}, Lemmas \ref{dilated-symbol-bound} and \ref{discretization} ensure that it is sufficient to prove that
\begin{equation}\label{aim-dis-Mg}
	\int_{\Sph}\left(\sup_{a_k\in\Gamma(\xi)}|g(a_k)|\right)
	^{\frac{pq}{p-q}}\dd\sigma(\xi)
	\less\|M_g\|_{\BT^p\to\BT^q}^{\frac{pq}{p-q}}
\end{equation}
when $Z=\{a_k\}\subset\B\setminus\{0\}$ is a $\delta$-lattice and $\delta$ is small enough. Assume $Z=\{a_k\}$ is such a lattice. For any $\lambda=\{\lambda_k\}\in T^p_{\infty}(Z)$, consider the test functions
\begin{equation*}
	S_{\tau}(\lambda)(z)
	:=\sum_{k=1}^{\infty}\lambda_kr_k(\tau)f_k(z),
\end{equation*}
where $\{r_k\}$ is a sequence of Rademacher functions on $[0,1]$, and
$$f_k(z)=\left(\frac{1-|a_k|^2}{1-\ip{z}{a_k}}\right)^{2\theta}-
\left(\frac{1-|a_k|^2}{1-\ip{z}{a_k}}\right)^{\theta}$$
with $\theta>n+n/p$. By Lemma \ref{synthesis}, $S_{\tau}(\lambda)\in\BT^p$ for almost every $\tau\in[0,1]$, and $\|S_{\tau}(\lambda)\|_{\BT^p}\less\|\lambda\|_{T^p_{\infty}(Z)}$. Note that for each $k\geq1$,
$$f_k(a_k)=0\quad \text{and}\quad Rf_k(a_k)=\frac{\theta|a_k|^2}{1-|a_k|^2}.$$
As in the proof of Theorem \ref{I_gp>q}, the boundedness of $M_g:\BT^p\to\BT^q$ then implies
\begin{align*}
\|M_g\|_{\BT^p\to\BT^q}^q\|\lambda\|_{T^p_{\infty}}^q
&\gtrsim\int_{\Sph}\int_0^1\sup_{z\in\Gamma(\xi)}
  \left|\sum_{k=1}^{\infty}\lambda_kr_k(\tau)
  R(gf_k)(z)(1-|z|^2)\right|^q\dd\tau\dd\sigma(\xi)\\
&\gtrsim\int_{\Sph}\left(\sup_{a_j\in\Gamma(\xi)}
  |\lambda_j||R(gf_j)(a_j)|(1-|a_j|^2)\right)^q\dd\sigma(\xi)\\
&\asymp\int_{\Sph}\left(\sup_{a_j\in\Gamma(\xi)}
  |\lambda_j||a_j|^2|g(a_j)|\right)^q\dd\sigma(\xi).
\end{align*}
Due to Lemmas \ref{discrete-tent-multipliers} and \ref{necessary}, this yields the desired estimate \eqref{aim-dis-Mg}.

(3) If $M_g$ is bounded, (2) gives $g\in H^{\frac{pq}{p-q}}$ with
$\|g\|_{H^{pq/(p-q)}}\lesssim\|M_g\|_{\BT^p\to\BT^q}$. Then by Theorem \ref{I_gp>q}, $I_g$ is bounded with $\|I_g\|_{\BT^p\to\BT^q}\lesssim\|g\|_{H^{pq/(p-q)}}$. Consequently, the identity \eqref{eq:identities} gives the boundedness of $J_g$ together with the estimate $\|J_g\|_{\BT^p\to\BT^q}\lesssim\|M_g\|_{\BT^p\to\BT^q}$. Conversely, if $g\in H^{\frac{pq}{p-q}}$ and $J_g$ is bounded, then Theorem \ref{I_gp>q} in conjunction with \eqref{eq:identities} gives
\[
 \|M_g\|_{\BT^p\to\BT^q}
 \lesssim\|J_g\|_{\BT^p\to\BT^q}+\|I_g\|_{\BT^p\to\BT^q}+|g(0)|
 \lesssim\|g\|_{H^{pq/(p-q)}}+\|J_g\|_{\BT^p\to\BT^q}.
\]
The proof is complete.
\end{proof}

By combining Theorems \ref{I_gp>q}, \ref{J_gp>q} and \ref{M_gp>q}, we complete the proof of Theorem \ref{main3}.

\section{Examples}\label{example}

In this section, we provide some counterexamples concerning the conditions in our main theorems. We will work on the unit disk $\D=\mathbb{B}_1$. In this setting, the definition of $J_g$ reduces to
$$J_gf(z)=\int_0^{z}f(w)g'(w)\dd w,\quad z\in\D,\quad f\in\hol(\D).$$
Throughout this section, we use the equivalent (quasi-)norms of the corresponding spaces defined with the ordinary derivative and use $\Q_{\gamma}$ accordingly.

The following example indicates that, at least in the unit disk setting, the necessary conditions in Theorem \ref{main2} (2) and Theorem \ref{main3} (2) are not sufficient.

\begin{example}\label{ex:bloch}
Take $\lambda_k=2^{2^k}$ and define
$$g(z)=\sum_{k=1}^{\infty}z^{\lambda_k},\quad z\in\D.$$
Then
$$g\in\bloch(\D)\subset\bigcap_{0<s<\infty}\BT^s(\D),$$
and $J_g:\BT^p(\D)\to\BT^q(\D)$ is not bounded for any $0<p,q<\infty$.
\end{example}
\begin{proof}
For each positive integer $K$, let $f_K(z)=\sum_{k=1}^K z^{\lambda_k}$. Note that for $z\in\D\setminus\{0\}$,
\begin{align*}
\sum_{k=1}^{\infty}\lambda_{k}|z|^{\lambda_k-1}
&=\sum_{k=1}^{\infty}\int_{\lambda_k}^{2\lambda_k}
	|z|^{\lambda_k-1}\dd t\\
&\leq\sum_{k=1}^{\infty}\int_{\lambda_k}^{2\lambda_k}
	|z|^{t/2-1}\dd t\\
&\leq\sum_{k=1}^{\infty}\int_{\lambda_k}^{\lambda_{k+1}}
	|z|^{t/2-1}\dd t\\
&=\int_{4}^{+\infty}|z|^{t/2-1}\dd t\\
&=-\frac{2|z|}{\log{|z|}}\less\frac{1}{1-|z|^2}.
\end{align*}
Thus $g\in\bloch(\D)\subset\bigcap_{0<s<\infty}\BT^s(\D)$ and $\sup_K\norm{f_K}_{\BT^p(\D)}<\infty$.
	
Put $r_K=1-\lambda_K^{-1}$. Then $r_K^{\lambda_K}\to e^{-1}$. For fixed $K\geq2$, write
$$S_K^-:=\sum_{k<K}\lambda_kr_K^{\lambda_k-1}
\quad\text{and}\quad
S_K^+:=\sum_{k>K}\lambda_kr_K^{\lambda_k-1}.$$
We claim that
\begin{equation}\label{SN}
\frac{S_K^-}{\lambda_K}\to0
\quad\text{and}\quad
\frac{S_K^+}{\lambda_K}\to0
\end{equation}
as $K\to\infty$. To see this, note that $\lambda_{k+1}=\lambda_k^2$. Consequently,
$$S_K^-\leq\sum_{k<K}\lambda_k\leq(K-1)\lambda_{K-1}=(K-1)\sqrt{\lambda_K},$$
which implies the first limit in \eqref{SN}. For the second limit, noting that
$$1-x\leq e^{-x},\quad 0<x<1,$$
we have
\begin{align*}
S_K^+
&=\sum_{k=1}^{\infty}\lambda_{K+k}
  (1-\lambda_K^{-1})^{\lambda_{K+k}-1}\\
&\less\sum_{k=1}^{\infty}\lambda_{K+k}
  \exp\left(-\frac{\lambda_{K+k}}{\lambda_K}\right)\\
&=\sum_{k=1}^{\infty}\lambda_K^{2^k}
  \exp\left(-\lambda_K^{2^k-1}\right)\\
&=\lambda_K^2e^{-\lambda_K}\left(1+
  \sum_{k=2}^{\infty}\lambda_K^{2^k-2}
  \exp(\lambda_K-\lambda_K^{2^k-1})\right)\\
&\leq\lambda_K^2e^{-\lambda_K}\left(1+
  \sum_{k=2}^{\infty}
  \exp(\lambda_K^{2^{k-1}-1}+\lambda_K-\lambda_K^{2^k-1})\right).
\end{align*}
It is clear that for $k\geq2$,
$$\lambda_K^{2^{k-1}-1}\leq\sqrt{\lambda_K^{2^k-1}}
\quad\text{and}\quad
\lambda_K\leq\sqrt{\lambda_K^{2^k-1}},$$
which gives
$$\lambda_K^{2^k-1}-(\lambda_K^{2^{k-1}-1}+\lambda_K)
\geq\lambda_K^{2^k-1}-2\sqrt{\lambda_K^{2^k-1}}
>\lambda_K^{2^{k-1}-1/2}$$
since $\lambda_K^{2^k-1}>9$. Consequently,
\begin{equation}\label{Sk+}
S_K^+\less\lambda_K^2e^{-\lambda_K}\left(1+
\sum_{k=2}^{\infty}\exp(-\lambda_K^{2^{k-1}-1/2})\right)
\less\lambda_K^2e^{-\lambda_K},
\end{equation}
which gives the second limit in \eqref{SN}. Therefore, noting that $r_K^{\lambda_K}\to e^{-1}$, we can choose $K$ so large that
$$S_K^-+S_K^+<\frac{1}{2}\lambda_Kr_K^{\lambda_K-1}.$$
This yields that for any $\xi\in\T$,
$$|g'(r_K\xi)|\geq\lambda_Kr_K^{\lambda_K-1}-
\left|\sum_{k\neq K}\lambda_k(r_K\xi)^{\lambda_k-1}\right|
\geq\lambda_Kr_K^{\lambda_K-1}-(S_K^-+S_K^+)
>\frac{1}{2}\lambda_Kr_K^{\lambda_K-1}.$$
Combining this with the fact that $r_K\xi\in\Gamma(\xi)$, we arrive at
\begin{align*}
\|J_gf_K\|^q_{\BT^q(\D)}
&=\int_{\T}\left(\sup_{z\in\Gamma(\xi)}
  |g'(z)f_K(z)|(1-|z|^2)\right)^q\dd\sigma(\xi)\\
&\geq\int_{\T}|g'(r_K\xi)f_K(r_K\xi)|^q(1-r_K^2)^q\dd\sigma(\xi)\\
&\gtrsim\int_{\T}|f_K(r_K\xi)|^q\dd\sigma(\xi).
\end{align*}
The classical lacunary series estimate (see \cite[Theorem 6.2.2]{JVA}) then yields
$$\|J_gf_K\|_{\BT^q(\D)}\gtrsim
\left(\sum_{k=1}^Kr_K^{2\lambda_k}\right)^{1/2}\asymp\sqrt{K}.$$
Thus $J_g:\BT^p(\D)\to\BT^q(\D)$ is not bounded.
\end{proof}

The following example indicates that the necessary conditions in Theorem \ref{main2} (4) and in Theorem \ref{main3} (4) are not sufficient in general.

\begin{example}\label{ex:blaschke}
Let $2<q<\infty$. There exists a Blaschke product $B$ on $\D$ that fails to define a bounded multiplier from $\BT^p(\D)$ to $\BT^q(\D)$ for all $0<p<\infty$.
\end{example}
\begin{proof}
	Let $f(z)=\sum_{k=1}^{\infty}z^{2^k}$. Then for any $z\in\D\setminus\{0\}$,
	$$|f'(z)|\leq\sum_{k=1}^{\infty}2^k|z|^{2^k-1}
	\leq\sum_{k=1}^{\infty}\int_{2^k}^{2^{k+1}}|z|^{t/2-1}\dd t
	=\int_2^{+\infty}|z|^{t/2-1}\dd t
	=-\frac{2}{\log|z|}\less\frac{1}{1-|z|^2}.$$
	That is, $f\in\bloch(\D)\subset\BT^p(\D)$ for any $0<p<\infty$.
	
	Fix $1<s<q/2$ and put
	$$r_k=1-2^{-k},\quad m_k=\left\lfloor\frac{2^k}{k^s}\right\rfloor.$$
	Choose $k_0$ sufficiently large so that $m_k\ge2$ for every $k\ge k_0$. For each $k\geq k_0$, choose a rotation parameter $\theta_k$, and set
	\[
	a_{k,j}=r_k\exp\left(i\theta_k+\frac{2\pi i j}{m_k}\right),
	\qquad 0\leq j<m_k.
	\]
	Averaging in $\theta_k$ allows a choice such that
	\begin{equation}\label{theta}
		\sum_{j=0}^{m_k-1}|f(a_{k,j})|^q
		\geq m_k\int_{\T}|f(r_k\xi)|^q\dd\sigma(\xi).
	\end{equation}
	Since $q>2$, H\"older's inequality gives that
	\begin{equation}\label{>k^p}
		\int_{\T}|f(r_k\xi)|^q\dd\sigma(\xi)
		\geq\left(\int_{\T}|f(r_k\xi)|^2\dd\sigma(\xi)\right)^{q/2}
		\geq\left(\sum_{l=1}^kr_k^{2^{l+1}}\right)^{q/2}
		\asymp k^{q/2}.
	\end{equation}
	It is clear that
	$$\sum_{k,j}(1-|a_{k,j}|)=\sum_{k=k_0}^{\infty}2^{-k}m_k
	\less\sum_{k=1}^{\infty}k^{-s}<\infty$$
	since $s>1$. Let $B$ be the Blaschke product with zero set $\mathcal{Z}=\{a_{k,j}\}_{k,j}$. We claim that the sequence $\mathcal{Z}$ is an interpolating sequence for $H^{\infty}(\D)$, whose proof is postponed. Then by the Carleson interpolation theorem (see, for instance, \cite[Theorem 1.1, p. 278]{Ga}), the Blaschke product $B$ satisfies
	\[
	(1-|a_{k,j}|^2)|B'(a_{k,j})|
	=\prod_{(\ell,j')\ne(k,j)}
	\rho(a_{k,j},a_{\ell,j'})\geq\delta
	\]
	for some $\delta>0$ and any $k,j$, and hence
	$$|B'(a_{k,j})|\geq\frac{\delta}{1-|a_{k,j}|^2},$$
	where
	\[
	\rho(z,w)=\left|\frac{z-w}{1-\overline wz}\right|
	\]
	is the pseudo-hyperbolic metric on $\D$. 
	In view of $(Bf)'(a_{k,j})=B'(a_{k,j})f(a_{k,j})$, we have
	\begin{equation}\label{Bf'}
		|(Bf)'(a_{k,j})|(1-|a_{k,j}|^2)\geq\delta|f(a_{k,j})|.
	\end{equation}
	Since $\sigma(I(a_{k,j}))\asymp1-|a_{k,j}|=2^{-k}$, we know that the intervals $I(a_{k,0})$, \dots, $I(a_{k,m_k-1})$ are pairwise disjoint whenever $k$ is large enough. This together with \eqref{theta}, \eqref{>k^p} and \eqref{Bf'} yields 
	\begin{align*}
		&\int_{\T}\left(\sup_{z\in\Gamma(\xi)}
		|(Bf)'(z)|(1-|z|^2)\right)^q\dd\sigma(\xi)\\
		&\ \geq\sum_{j=0}^{m_k-1}\int_{I(a_{k,j})}
		\left(\sup_{z\in\Gamma(\xi)}
		|(Bf)'(z)|(1-|z|^2)\right)^q\dd\sigma(\xi)\\
		&\ \geq\sum_{j=0}^{m_k-1}\int_{I(a_{k,j})}
		\left(|(Bf)'(a_{k,j})|(1-|a_{k,j}|^2)\right)^q\dd\sigma(\xi)\\
		&\ \gtrsim\sum_{j=0}^{m_k-1}2^{-k}|f(a_{k,j})|^q\\
		&\ \gtrsim 2^{-k}m_kk^{q/2}\asymp k^{q/2-s}\to\infty
	\end{align*}
	as $k\to\infty$ since $s<q/2$. Therefore, $Bf\notin\BT^q(\D)$.
	
	It remains to show that the sequence $\mathcal{Z}$ is an interpolating sequence for $H^{\infty}(\D)$. To this end, we first verify uniform separation. If
	$k<\ell$, then $r_k<r_\ell$ and hence
	\[
	\rho(a_{k,j_1},a_{\ell,j_2})
	\geq\frac{r_\ell-r_k}{1-r_kr_\ell}
	=\frac{2^{-k}-2^{-\ell}}
	{2^{-k}+2^{-\ell}-2^{-k-\ell}}
	\geq\frac13.
	\]
	For two distinct points $a_{k,j_1}$ and $a_{k,j_2}$ in the same $k$-th layer, put
	$$\alpha=\frac{2\pi}{m_k}\min\{|j_1-j_2|,m_k-|j_1-j_2|\}.$$
	Then $\alpha\in[2\pi/m_k,\pi]$. A direct computation gives
	\[
	\rho(a_{k,j_1},a_{k,j_2})^2
	=\frac{4r_k^2\sin^2(\alpha/2)}
	{(1-r_k^2)^2+4r_k^2\sin^2(\alpha/2)}.
	\]
	For all $k\geq k_0$, $r_k\geq1/2$ and
	$\sin(\alpha/2)\geq\sin(\pi/m_k)\geq2/m_k$. Since
	$1-r_k^2\leq2\cdot2^{-k}$ and
	$m_k2^{-k}\leq k^{-s}$, it follows that
	\[
	\frac{2r_k\sin(\alpha/2)}{1-r_k^2}
	\geq\frac{1}{m_k2^{-k}}\geq1.
	\]
	Thus $\rho(a_{k,j_1},a_{k,j_2})\geq2^{-1/2}$ for all $k\geq k_0$. Consequently, we obtain
	\begin{equation}\label{uniform-sep}
		\inf_{(k,j_1)\ne(\ell,j_2)}\rho(a_{k,j_1},a_{\ell,j_2})\geq\frac{1}{3}.
	\end{equation}
	We next verify the Carleson measure condition. Let $I\subset\T$ be an arc
	of normalized length $|I|$, and let
	\[
	S(I)=\{re^{it}:e^{it}\in I,\ 1-|I|\leq r<1\}
	\]
	be its Carleson square. If $a_{k,j}\in S(I)$, then $2^{-k}=1-r_k\leq |I|$.
	Moreover, because the $m_k$ points on the $k$-th circle are equally spaced,
	\[
	\#\{j:a_{k,j}\in S(I)\}\leq m_k|I|+1.
	\]
	Therefore,
	\begin{align*}
		\sum_{a_{k,j}\in S(I)}(1-|a_{k,j}|^2)
		&\leq2\sum_{2^{-k}\leq |I|}2^{-k}(m_k|I|+1)\\
		&\leq2|I|\sum_{k=k_0}^{\infty}m_k2^{-k}
		+2\sum_{2^{-k}\leq |I|}2^{-k}\\
		&\lesssim |I|\sum_{k=1}^{\infty}k^{-s}+|I|
		\lesssim |I|.
	\end{align*}
	Thus
	\[
	\mu=\sum_{k,j}(1-|a_{k,j}|^2)\delta_{a_{k,j}}
	\]
	is a Carleson measure. By \cite[Theorem 1.1, p. 278]{Ga}, the pseudo-hyperbolic separation in \eqref{uniform-sep} together with this Carleson
	condition implies that $\mathcal Z$ is an interpolating sequence for
	$H^\infty(\D)$, which finishes the proof.
\end{proof}

The following example indicates that, at least in the unit disk setting, the sufficient conditions in Theorem \ref{main3} (3) and (5) are not necessary. Since $\bloch_{\log}(\D)\subset\bigcap_{0<t<\infty}\BT^t_{\log}(\D)$, this also shows that the sufficient conditions in Theorem \ref{main2} (3) and (5) are not necessary.

\begin{example}\label{lacunary-multiplier}
	Let $\lambda_k=2^{2^k}$, $c_k=2^{-3k/4}$, and define
	\[
	g(z)=\sum_{k=1}^{\infty}c_kz^{\lambda_k},\quad z\in\D.
	\]
	Then $g\in H^\infty(\D)\setminus\bigcup_{0<t<\infty}\BT^t_{\log}(\D)$, whereas both $M_g$ and $J_g$ are bounded from $\BT^p(\D)$ to $\BT^q(\D)$ for every $0<q\leq p<\infty$.
\end{example}

To prove this, we need two lemmas. The first one estimates the truncated maximal functions on $\BT^p(\D)$, and the second one estimates the monomial multipliers on $\BT^p(\D)$.

\begin{lemma}\label{truncated-maximal}
	Let $0<p<\infty$ and $0<r<1$. For $f\in\BT^p(\D)$, put
	\[
	\N^{(r)}f(\xi):=\sup_{z\in\Gamma(\xi),\ |z|\leq r}|f(z)|.
	\]
	Then
	\begin{equation}\label{sqrt-log-estimate}
		\|\N^{(r)}f\|_{L^p(\T)}\lesssim
		\left(\log\frac 2{1-r}\right)^{1/2}\|f\|_{\BT^p(\D)}
	\end{equation}
	with implicit constant independent of $r$.
\end{lemma}
\begin{proof}
For $\gamma>1$ and $\xi\in\T$, write
\[
\Q_{\gamma}f(\xi):=\sup_{z\in\Gamma_{\gamma}(\xi)}|f'(z)|(1-|z|^2)
\]
and
$$\N_{\gamma}f(\xi):=\sup_{z\in\Gamma_{\gamma}(\xi)}|f(z)|.$$
We first record the change-of-aperture estimate
\begin{equation}\label{aperture-maximal}
\|\Q_{\gamma}f\|_{L^p(\T)}\lesssim\|\Q f\|_{L^p(\T)}.
\end{equation}
For $r\leq1/2$, \eqref{sqrt-log-estimate} follows from the local growth estimate \eqref{growth}. Suppose $r>1/2$, set $\rho=(1+r)/2$, and write $f_{\rho}(z)=f(\rho z)$. Define the analytic area function by
\[
\mathcal S(f_{\rho})(\xi):=
\left(\int_{\Gamma(\xi)}|(f_{\rho})'(w)|^2\dd A(w)\right)^{1/2},
\quad \xi\in\T,
\]
where $dA$ is the normalized area measure on $\D$. If $w\in\Gamma(\xi)$, then $\rho w\in\Gamma_{4}(\xi)$, since
\[
|1-\rho w\overline\xi|
\leq(1-\rho)+\rho|1-w\overline\xi|
\leq2(1-\rho^2|w|^2).
\]
It follows that
\[
|(f_{\rho})'(w)|=|\rho f'(\rho w)|
\leq\frac{\Q_{4}f(\xi)}{1-\rho^2|w|^2}.
\]
Note that the angular width of $\Gamma(\xi)$ at radius $t$ is at most $C(1-t)$. Hence
\begin{align*}
\int_{\Gamma(\xi)}|(f_{\rho})'(w)|^2dA(w)
&\less\int_{0}^1\int_{\{\zeta\in\T:t\zeta\in\Gamma(\xi)\}}
	|(f_{\rho})'(t\zeta)|^2\dd\sigma(\zeta)\dd t\\
&\lesssim \Q_{4}f(\xi)^2
	\int_{0}^1\frac{1-t}{(1-\rho^2t^2)^2}\dd t\\
&\leq \Q_{4}f(\xi)^2
	\int_{0}^1\frac{dt}{1-\rho^2t^2}\\
&\lesssim \Q_{4}f(\xi)^2\log\frac{2}{1-\rho}.
\end{align*}
Consequently,
$$\mathcal{S}(f_{\rho})(\xi)\lesssim
\Q_{4}f(\xi)\left(\log\frac{2}{1-\rho}\right)^{1/2}.$$
The classical Calder\'on area theorem (see \cite{Ca,FS} or \cite[Theorem B]{Pa}) together with \eqref{growth} and \eqref{aperture-maximal} therefore implies
\begin{equation}\label{dilated-Hardy-estimate}
\|f_\rho\|_{H^p(\D)}\asymp
|f_{\rho}(0)|+\|\mathcal{S}(f_{\rho})\|_{L^p(\T)}
\lesssim
\left(\log\frac 2{1-\rho}\right)^{1/2}\|f\|_{\BT^p(\D)}.
\end{equation}
Suppose now that $z\in\Gamma(\xi)$ and $|z|\leq r$. Since $0<r<1$ and $|z|\leq r$, we have
\[
1-\frac{|z|}{\rho}=1-\frac{2}{1+r}|z|\geq\frac{1-|z|}{2},
\]
which implies
\[
\left|1-\frac{z}{\rho}\overline\xi\right|
\leq |1-z\overline\xi|+\left(\frac{1}{\rho}-1\right)|z|
\leq2(1-|z|)+\frac{1-r}{1+r}
<3(1-|z|)\leq6(1-|z/\rho|^2).
\]
That is, $z/\rho\in\Gamma_{12}(\xi)$, and so
$$\N^{(r)}f\leq \N_{12}(f_\rho).$$
Applying the standard non-tangential maximal estimate for Hardy spaces (see, for instance, \cite[Theorem 4.24]{Zh}) and \eqref{dilated-Hardy-estimate}, and using the equality $1-\rho=(1-r)/2$, we establish
$$\|\N^{(r)}f\|_{L^p(\T)}\leq\|\N_{12}(f_{\rho})\|_{L^p(\T)}
\lesssim\|f_{\rho}\|_{H^p(\D)}
\lesssim\left(\log\frac 2{1-r}\right)^{1/2}
\|f\|_{\BT^p(\D)},$$
and the proof is complete.
\end{proof}

\begin{lemma}\label{monomial-multipliers}
	Let $0<p<\infty$. For every integer $k\geq1$,
	\[
	\|M_{z^k}\|_{\BT^p(\D)\to\BT^p(\D)}
	\lesssim\sqrt{\log(ek)},
	\]
	where the implicit constant depends only on $p$.
\end{lemma}
\begin{proof}
	Since $(J_{z^k}f)'(z)=kz^{k-1}f(z)$,
	\[
	\Q({J_{z^k}f})(\xi)=\sup_{z\in\Gamma(\xi)}
	k(1-|z|^2)|z|^{k-1}|f(z)|.
	\]
	Decompose the disk into the annuli
	\[
	A_j:=\{z:2^{-j-1}<1-|z|\leq2^{-j}\},\qquad j\geq0,
	\]
	and put $r_j=1-2^{-j-1}$. For $z\in A_j$,
	\[
	k(1-|z|^2)|z|^{k-1}\leq2k2^{-j}(1-2^{-j-1})^{k-1}
	\leq2k2^{-j}\exp(-(k-1)2^{-j-1}).
	\]
	Therefore,
	\[
	\Q({J_{z^k}f})(\xi)\leq
	2\sum_{j\geq0}k2^{-j}e^{-(k-1)2^{-j-1}}\N^{(r_j)}f(\xi).
	\]
	Put $s=\min\{1,p\}$ and $w_{k,j}=k2^{-j}e^{-(k-1)2^{-j-1}}$.
	Since $J_{z^k}f(0)=0$, Minkowski's inequality when $p\geq1$, and
	$p$-subadditivity when $p<1$, together with Lemma~\ref{truncated-maximal}, give
	\[
	 \|J_{z^k}f\|_{\BT^p(\D)}^s
	 \leq2^s\sum_{j\geq0}w_{k,j}^s\|\N^{(r_j)}f\|_{L^p(\T)}^s
	 \lesssim\|f\|_{\BT^p(\D)}^s
	 \sum_{j\geq0}w_{k,j}^s(j+1)^{s/2}.
	\]
	Set $L=\lceil\log_2k\rceil$. For $j=L+m\geq L$, $w_{k,j}\leq2^{-m}$,
	whereas for $j=L-m<L$, $w_{k,j}\leq2^me^{-2^{m-3}}$. Consequently,
	$$\sum_{j\geq L}w_{k,j}^s(j+1)^{s/2}
	\leq(L+1)^{s/2}\sum_{m\geq0}2^{-ms}(m+1)^{s/2}
	\lesssim(L+1)^{s/2},$$
	and
	$$\sum_{j<L}w_{k,j}^s(j+1)^{s/2}
	\leq(L+1)^{s/2}\sum_{m\geq1}2^{ms}e^{-s2^{m-3}}
	\lesssim(L+1)^{s/2}.$$
	Hence
	\[
	\|J_{z^k}f\|_{\BT^p(\D)}\lesssim\sqrt{L+1}\|f\|_{\BT^p(\D)}
	\lesssim\sqrt{\log(ek)}\|f\|_{\BT^p(\D)}.
	\]
	Since the proof of Theorem \ref{main2} gives $\|I_{z^k}\|_{\BT^p(\D)\to\BT^p(\D)}\leq1$, the equality
	$M_{z^k}=I_{z^k}+J_{z^k}$ yields the desired multiplier estimate.
\end{proof}

We are now ready to prove Example \ref{lacunary-multiplier}.

\begin{proof}[Proof of Example \ref{lacunary-multiplier}]
	Since $\sum_{k=1}^{\infty}c_k<\infty$, the series converges uniformly on the closed disk $\overline{\D}$ and defines a function in $H^\infty(\D)$.
	Fix $0<p<\infty$ and put $s=\min\{1,p\}$. Then
	\[
	 \sum_{k\geq1}c_k^s[\log(e\lambda_k)]^{s/2}
	 \lesssim\sum_{k\geq1}2^{-3ks/4}2^{ks/2}
	 =\sum_{k\geq1}2^{-ks/4}<\infty.
	\]
	For $g_K(z):=\sum_{k=1}^Kc_kz^{\lambda_k}$, note that $g_K(0)=0$. Minkowski's inequality if $p\geq1$, or $p$-subadditivity if $p<1$,
	and Lemma~\ref{monomial-multipliers} yield, for $f\in\BT^p(\D)$,
	\[
	 \|\Q(g_Kf)\|^s_{L^p(\T)}=\|g_Kf\|_{\BT^p(\D)}^s
	 \leq\sum_{k=1}^K c_k^s\|z^{\lambda_k}f\|_{\BT^p(\D)}^s
	 \lesssim\left(\sum_{k\geq1}2^{-ks/4}\right)\|f\|_{\BT^p(\D)}^s.
	\]
	The functions $g_Kf$ and their derivatives converge locally uniformly to
	$gf$ and its derivative. Consequently,
	\[
	 \Q({gf})\leq\liminf_{K\to\infty}\Q({g_Kf}).
	\]
	Fatou's lemma together with the fact $g(0)=0$ yields
	\[
	 \|M_gf\|_{\BT^p(\D)}=\|\Q({gf})\|_{L^p(\T)}
	 \lesssim\left(\sum_{k\geq1}2^{-ks/4}\right)^{1/s}\|f\|_{\BT^p(\D)}.
	\]
	Thus $M_g$ is bounded on $\BT^p(\D)$ for every $p>0$.
	Since $g\in H^\infty$, Theorem~\ref{main2} and
	\eqref{eq:identities} also give the boundedness of $J_g$ on $\BT^p(\D)$. For $q<p$, the inclusion $\BT^p(\D)\subset\BT^q(\D)$ is continuous. Hence both $M_g:\BT^p(\D)\to\BT^q(\D)$ and $J_g:\BT^p(\D)\to\BT^q(\D)$ are bounded as well.

We next prove that
\begin{equation}\label{infty-everywhere}
\sup_{z\in\Gamma(\xi)}|g'(z)|(1-|z|^2)\log\frac2{1-|z|^2}
=\infty,\qquad\forall\xi\in\T,
\end{equation}
which implies $g\notin\bigcup_{0<t<\infty}\BT^t_{\log}(\D)$ and will finish the proof. Put $r_k=1-\lambda_k^{-1}$. For each $K\geq2$, it is clear that
$$S_K^-:=\sum_{k<K}c_k\lambda_kr_K^{\lambda_k-1}
\leq\lambda_{K-1}\sum_{k<K}c_k\lesssim\sqrt{\lambda_{K}},$$
which implies
$$\frac{S_K^-}{c_K\lambda_Kr_K^{\lambda_K-1}}
\less\frac{2^{3K/4}}{\sqrt{\lambda_{K}}}\to0$$
as $K\to\infty$. For $k>K$, $c_k\leq c_K$, and so by \eqref{Sk+},
$$S_K^+:=\sum_{k>K}c_k\lambda_kr_K^{\lambda_k-1}
\less c_K\lambda_K^2e^{-\lambda_K},$$
which gives
$$\frac{S_K^+}{c_K\lambda_Kr_K^{\lambda_K-1}}
\less\frac{\lambda_Ke^{-\lambda_K}}{r_K^{\lambda_K-1}}\to0$$
as $K\to\infty$. Hence we can choose $K$ so large that
$$S_K^-+S_K^+<\frac{1}{2}c_K\lambda_Kr_K^{\lambda_K-1}.$$
Then for any $\xi\in\T$, the triangle inequality gives
$$|g'(r_K\xi)|\geq c_K\lambda_Kr_K^{\lambda_K-1}-S_K^--S_K^+
\gtrsim c_K\lambda_K.$$
Since $(1-r_K^2)\lambda_K\asymp1$ and
$\log(2/(1-r_K^2))\asymp\log\lambda_K$, it follows uniformly in $\xi\in\T$ that
\[
|g'(r_K\xi)|(1-r_K^2)\log\frac2{1-r_K^2}
\gtrsim c_K\log\lambda_K=2^{K/4}\log2\to\infty
\]
as $K\to\infty$. Since $r_k\xi\in\Gamma(\xi)$ for each $k\geq1$, we establish \eqref{infty-everywhere} and finish the proof.
\end{proof}

We conclude this section with an example showing that $\bloch_{\log}(\D)\subsetneq\bigcap_{0<p<\infty}\BT^p_{\log}(\D)$.

\begin{example}\label{ex:logarithmic-symbol}
Let
\[
 g(z)=\log\frac1{1-z},\qquad z\in\D,
\]
where we use the holomorphic logarithm on the right half-plane with $\log1=0$. Then
\[
 g\notin H^\infty(\D)\cup\bloch_{\log}(\D),
 \qquad g\in\bigcap_{0<p<\infty}\BT_{\log}^p(\D).
\]
\end{example}
\begin{proof}
Since $g(t)=\log(1/(1-t))\to\infty$ as $t\to1^-$, the function $g$
is not bounded. Moreover,
\[
 |g'(t)|(1-t^2)\log\frac2{1-t^2}=(1+t)\log\frac{2}{1-t^2}\to\infty
\]
as $t\to1^-$, and hence $g\notin\bloch_{\log}(\D)$.

For $z\in\D$, we have
$1-|z|^2\leq2|1-z|$ and $0<|1-z|<2$. Set $x=(1-|z|^2)/(2|1-z|)\in(0,1]$.
Since $x\log(1/x)\leq1/e$, we obtain
\begin{align*}
|g'(z)|(1-|z|^2)\log\frac{2}{1-|z|^2}
&=\frac{1-|z|^2}{|1-z|}\log\frac{2}{1-|z|^2}\\
&=2x\log\frac1{x|1-z|}\\
&=2x\log\frac{4}{|1-z|}+2x\log\frac1x-2x\log4\\
&\leq2\log\frac{4}{|1-z|}+\frac2e\\
&\lesssim\log\frac{4}{|1-z|}.
\end{align*}
For $\xi\in\T$, if $z\in\Gamma(\xi)$, Lemma \ref{kernel-estimate}, applied with
$a=t\in(0,1)$ and followed by $t\to1^-$, gives
\[
 |1-\xi|\leq9|1-z|.
\]
Consequently, for $\xi\in\T\setminus\{1\}$,
\begin{equation}\label{log-symbol-majorant}
 \sup_{z\in\Gamma(\xi)}
 |g'(z)|(1-|z|^2)\log\frac2{1-|z|^2}
 \lesssim\log\frac{36}{|1-\xi|}.
\end{equation}
For any $0<p<\infty$,
\begin{align*}
\int_{\T}\left(\log\frac{36}{|1-\xi|}\right)^p\dd\sigma(\xi)
&\less1+\int_{-1}^{1}
  \left(\log\frac{36}{|1-e^{i\vartheta}|}\right)^p\dd\vartheta\\
&\less1+
  \int_{-1}^1\left(\log\frac{2}{|\vartheta|}\right)^p\dd\vartheta\\
&\less1+\int_0^{+\infty}t^pe^{-t}\dd t<\infty.
\end{align*}
This together with \eqref{log-symbol-majorant} proves $g\in\BT_{\log}^p(\D)$ for every $0<p<\infty$.
\end{proof}

\noindent{\bf AI Assistance Disclosure.} OpenAI's GPT-6 Astra provided mathematical suggestions that contributed to the development of the proofs of Lemmas \ref{synthesis} and \ref{dilated-symbol-bound}. The explicit construction in Section \ref{example} were also developed with the assistance of GPT-6 Astra. All AI-assisted contributions have been carefully reviewed and verified by the author, who takes full responsibility for the content of the manuscript.


\medskip 



\end{document}